\documentclass[final,1p,times,authoryear]{elsarticle}
\usepackage{amsmath}
\usepackage{amsthm}
\usepackage{amssymb}
\usepackage{amsfonts}
\usepackage{algorithmicx,algpseudocode}
\usepackage[section]{algorithm}

\newtheorem{theorem}{Theorem}

\theoremstyle{definition} \newtheorem{definition}{Definition}[section]

\newcommand{\Hom}{{\rm Hom}}

\usepackage{amssymb}
\usepackage{graphicx}
\usepackage{amsfonts}
\usepackage{algorithmicx,algpseudocode}
\usepackage[section]{algorithm}
\theoremstyle{definition} \newtheorem{datatype}{Data type}[section]

\begin{document}
\begin{frontmatter}

\title{Cohomology with local coefficients
 and knotted manifolds}

\author{Graham Ellis}
\address{School of Mathematics, National University of Ireland, Galway}
\ead{graham.ellis@nuigalway.ie}
\ead[url]{http://hamilton.nuigalway.ie}

\author{Kelvin Killeen}
\address{School of Mathematics, National University of Ireland, Galway}
\ead{kelvin.killeen@nuigalway.ie}
\ead[url]{}

\begin{abstract}
We show how the classical notions of cohomology with local coefficients, 
CW-complex, covering space, homeomorphism equivalence, simple homotopy equivalence, 
tubular neighbourhood, and spinning can be encoded on a computer and used 
to calculate ambient isotopy invariants of continuous embeddings $N\hookrightarrow M$ of 
one topological manifold into another.
More specifically, we describe an algorithm for computing the homology $H_n(X,A)$ and
cohomology
$H^n(X,A)$ of a finite connected CW-complex X with local coefficients
in a $\mathbb Z\pi_1X$-module $A$ when $A$ is finitely generated
over $\mathbb Z$. It can be used, in particular, to compute the
integral cohomology $H^n(\widetilde X_H,\mathbb Z)$ and induced homomorphism
$H^n(X,\mathbb Z) \rightarrow H^n(\widetilde X_H,\mathbb Z)$ for the
covering map
$p\colon \widetilde X_H \rightarrow X$ associated to a finite index subgroup
$H < \pi_1X$, as well as the  corresponding homology homomorphism.
 We illustrate an  open-source implementation of the algorithm by using it to
show that: (i) the degree $2$  homology group $H_2(\widetilde X_H,\mathbb Z)$
distinguishes between the homotopy types of the complements  $X\subset \mathbb R^4$ of
 the spun Hopf link and Satoh's tube map of the welded Hopf link (these two complements having isomorphic fundamental groups and integral homology);
(ii) the degree $1$  homology homomorphism $H_1(p^{-1}(B),\mathbb Z)
\rightarrow H_1(\widetilde X_H,\mathbb Z)$ distinguishes between the
homeomorphism types of the complements $X\subset \mathbb R^3$ of the granny knot and the reef knot, where $B\subset X$ is the knot boundary (these two complements again having isomorphic fundamental groups and integral homology).
Our open source implementation allows the user to experiment with further examples of knots,  
knotted surfaces, and other embeddings of spaces.
 We conclude the paper with an explanation of how the cohomology algorithm
also  provides an approach  to computing the set
$[W,X]_\phi$ of based homotopy classes of maps
$f\colon W\rightarrow X$ of finite  CW-complexes over a fixed group homomorphism $\pi_1f=\phi$
in the
case where $\dim W =n$,  $\pi_1X$ is finite and $\pi_iX=0$ for $2\le i\le n-1$.

\end{abstract}

\end{frontmatter}

\section{Introduction}
Let $f,g\colon N\rightarrow M$ be two continuous cellular embeddings of a finite
CW-complex $N$ into a finite CW-complex $M$. A continuous cellular map 
$H\colon M\times [0,1] \rightarrow M, x\mapsto H_t(x)$ is said to be an {\em ambient isotopy} between $f$ and $g$ if $H_0$ is the identity map on $M$, each $H_t$ is a homeomorphism from $M$ to itself, and $H_1\circ f = g$.
 We are interested in computing invariants of embeddings with a view to
 distinguishing between their isotopy classes.  
The invariants are of most interest when $M$ and $N$ are manifolds. For embeddings of spheres $\mathbb S^n \rightarrow \mathbb S^{n+k}$ interest is on the case $k=2$ thanks to a result of \cite{MR123335} which states that  any two embeddings are ambient isotopic (in the piece-wise linear setting) for $k\ge 3$, and  a result of  
\cite{MR117693} which implies that there is again just one isotopy class of embeddings for $k=1$.  However, Zeeman also
 shows for instance that there is 
 more than one embedding $\mathbb S^n \rightarrow \mathbb S^1\times \mathbb S^{2n}$, up to ambient isotopy, for $n\ge 1$ (again in the piece-wise linear setting). 

The invariants we wish to calculate are the integral homology groups 
$H_k(\widetilde X_H, \mathbb Z)$ of finite covers of the complement 
$X=M\setminus N$, where
$H$ denotes the finite index subgroup of $\pi_1X$ arising as  
 the image of the induced homomorphism of fundamental groups 
$p_H\colon \pi_1\widetilde X_H \rightarrow \pi_1X$. By considering  
 a small open tubular  neighbourhood 
$N\subset N_\epsilon$ with boundary $\partial N_\epsilon$
we can consider the inclusion map $\partial N_\epsilon \hookrightarrow X$.
Letting $B =p_H^{-1}(\partial N_\epsilon)$ denote the preimage in $\widetilde {X_H} $ of the boundary, we also wish to calculate the induced homology homomorphism
$H_k(B,\mathbb Z) \rightarrow H_k(\widetilde X_H,\mathbb Z)$ as a means of distinguishing between embeddings of manifolds.

The use of homology of finite covers as an isotopy invariant is well documented in the literature. Indeed,   the first integral homology, $H_1(\widetilde{(\mathbb S^3\setminus \mathbb S^1)}_H,\mathbb Z)$, of finite coverings of knot complements
 was one of the earliest invariants used in the study of 
knot embeddings  $\mathbb S^1 \hookrightarrow \mathbb S^3$. Accounts of  methods for computing these first homology groups  directly 
from planar knot diagrams can be found, for instance, in  \cite{MR62125}, \cite{MR143201}, \cite{MR739142}, \cite{MR984809}.

The utility of  first homology of finite coverings in knot theory suggests an 
 analogous role for higher homology of coverings in the study of higher dimensional embeddings, especially those of codimension 2. Our aim is to describe various computer procedures for investigating the potential of this analogous role.  

Our procedures involve the
cellular chain complex $C_\ast\widetilde X$ of the universal cover
of a connected CW-complex $X$. This is a chain complex of free
$\mathbb ZG$-modules and $\mathbb ZG$-equivariant module homomorphisms
with $G=\pi_1X$ the fundamental group of $X$.
 Given a $\mathbb ZG$-module $A$, the  {\em homology} and {\em cohomology}
of $X$ {\em with local coefficients in} $A$ is defined as
\begin{equation}
H_n(X,A) = H_n(\,C_\ast \widetilde X \otimes_{\mathbb ZG} A\,) {\rm~~and~~}
H^n(X,A) = H^n(\,\Hom_{\mathbb ZG}(C_\ast \widetilde X , A)\,) .\label{EQcohom}
\end{equation}
Note that $H_n(\widetilde X_H,\mathbb Z) = H_n(X,A)$ when $A=\mathbb ZG\,\otimes_{\mathbb ZH} \mathbb Z$.
We present an algorithm for computing (\ref{EQcohom}) in the case where $A$ is a finitely generated $\mathbb ZG$-module. 
The algorithm requires a  CW-complex $X$ and details of $A$ as input.
Thus, in order to apply the algorithm to questions concerning embeddings of manifolds we need to provide algorithms for converting various classical
mathematical descriptions of
such embeddings to  embeddings of  CW-complexes. We provide algorithms for
this conversion for the case of  
 embeddings $N\hookrightarrow \mathbb S^{n+2}$ of an $n$-dimensional manifold $N$ 
into an $(n+2)$-sphere with $n=1$ and $n=2$, namely for the cases of links and surface links. 
  
The paper is organized as follows. In Section \ref{SECrep} we discuss the computer representation of CW-complexes. In Section \ref{SECnaive} we demonstrate some limitations of a naive approach to computing with (regular) CW-complexes and 
illustrate the need for implementations of homeomorphism equivalence and homotopy equivalence.
The  data types,  algorithms, and  mathematical theory underlying
the  computations in this section are explained in Sections \ref{SECchn} and \ref{SECdata}.
 In Section \ref{SECprojection} we consider the general principles involved in
 converting planar link diagrams to CW-structures on the complements of links in $\mathbb S^3$, and in converting surface link diagrams to CW-structures on the complements of surface links in $\mathbb S^4$.
 In Section \ref{SEChopf} we
provide an  illustrative computation
concerning  complements of two surfaces $K_1, K_2\subset \mathbb R^4$
 embedded in $\mathbb R^4$.  The surface $K_1$ is obtained by spinning the Hopf link in $\mathbb R^3$ about a plane  that does not intersect the link. The surface $K_2$ is obtained by applying the tube map of  \cite{MR1758871} to the welded Hopf link. The complements $X_i=\mathbb R^4\setminus K_i$ have isomorphic fundamental groups and isomorphic integral homology groups.
 It was shown by  \cite{MR2441256} that the spaces $X_1$, $X_2$ are homotopy inequivalent;
 their technique involves the fundamental crossed module derived from the lower dimensions of the universal cover of a space, and counts the representations of this fundamental crossed module into a given finite crossed module.
 We recover their homotopy inequivalence by using our algorithm to compute
$H_2(X_i , A)$ for $A=\mathbb Z\pi_1X_i \otimes_{\mathbb ZH}\mathbb Z$ with
$H<\pi_1X_i$ a certain finite index subgroup. We intentionally use an 
inefficient cubical complex representation of the spaces $X_1$, $X_2$ to demonstrate that algorithms are practical even for CW-complexes involving a large number of cells.  
In Section \ref{SECspin} we describe an algorithm that returns a regular CW-structure on the complement of the surface $S(L)$ in $\mathbb R^4$ obtained by
spinning a link $L$ about a plane which does not intersect $L$.
In Section \ref{SECgranny} we consider the granny knot $\kappa_1$ and reef knot $\kappa_2$.
The complements $X_i=\mathbb R^3\setminus \kappa_i$ are well known to have isomorphic fundamental groups and  different homeomorphism types. One way to distinguish between their homeomorphism types is to use the theory of quandles (see for instance \cite{MR3896314}).  We show that the homeomorphism types of these two complements can also be distinguished by the cokernels of the homology homomorphisms
$h_i\colon H_1(p_i^{-1}(B_i),\mathbb Z) \rightarrow H_1(\widetilde X_{i,H},\mathbb Z)$ where $p_i\colon \widetilde X_{i,H} \rightarrow X_i$ ranges over $6$-fold covering maps, and $B_i \subset X_i$ is the boundary of $\kappa_i$. As mentioned above, the  use of the pair $(X,B)$ and of the homology of finite covers are standard tools in knot theory. Nevertheless, it seems that this method of
distinguishing between the granny and reef knots is new.
 In Section \ref{SECnbhd} we describe an algorithm
that returns a regular CW-structure on
 the complement of an open tubular neighbourhood of a CW-subcomplex.

 The computation of cohomology with local coefficients has applications other than as an isotopy invariant.  In Section \ref{SECclass} we establish
 the following Theorem \ref{thm} which was the original
 motivation for this work.
 The theorem is a fairly immediate consequence of an old result of J.H.C.\,Whitehead and concerns the set $[W,X]_\phi$ of based homotopy classes of maps
$f\colon W\rightarrow X$ between connected
CW-complexes $W$ and $X$  that induce a
given homomorphism $\phi\colon \pi_1W\rightarrow \pi_1X$ of fundamental groups.
We assume that the spaces $W$, $X$  have preferred base points, that the
maps $f$ preserve base points, and that base points are preserved at each stage of  a homotopy between maps.
 The homotopy groups of $X$ are denoted by $\pi_nX$, $n\ge 0$.

\begin{theorem} \label{thm}
Let $W$ and $X$ be connected CW-complexes and let $\phi\colon \pi_1W
\rightarrow \pi_1X$ be a homomorphism between their fundamental groups.
Suppose that $\dim W=n$ and that
$\pi_iX=0$ for $2\le i\le n-1$. Then there is a non-canonical bijection
$$[W,X]_\phi \cong H^n(W, H_n(\widetilde X))$$
where the right hand side denotes cohomology with local coefficients in which
 $\pi_1W$ acts on  $H_n(\widetilde X)$
via the homomorphism $\phi$ and the canonical action of $\pi_1X$ on $H_n(\widetilde X)$.
\end{theorem}

We say that a CW-complex is {\em finite} if it contains only finitely many cells.
 Let us add to the hypothesis of Theorem \ref{thm} that $W$, $X$ are regular,
that $W$, $X$ and $\pi_1X$ are finite,  that $\phi\colon\pi_1W \rightarrow
\pi_1X$ is specified in terms of the free presentations of $\pi_1W$, $\pi_1X$
afforded by our computer implementation, and that some faithful
permutation representation
$\pi_1X \hookrightarrow S_k$ is at hand. Then
$H_n(\widetilde X)$ is finitely generated over $\mathbb Z$, and our \textsf{HAP}
implementation can in principle compute the cohomology group
$H^n(W, H_n(\widetilde X))$. The complexity of the computation is polynomial in the number of cells in $W$ and $X$. (This
 easy observation reduces to complexity bounds for Krushkal's spanning tree algorthm and for the Smith Normal Form algorithm.)
Hence one can  compute the set $[W,X]_\phi$ in polynomial time (ignoring practical constraints such as availability of memory).
This is a modest addition  to deeper results  in   \cite{CADEK:2014}, \cite{Krazal:2013}, \cite{kcral3}
 which establish the existence of a polynomial-time algorithm for computing
the set $[W,X]$ of homotopy classes of maps $W\rightarrow X$ under the
hypothesis that $W$, $X$ are simplicial sets, that $\pi_iX=0$ for $0\le i\le n-1$, and that $\dim W\le 2n-2$.

Section \ref{SECrepro} provides details on how to reproduce the computations of the paper.

\section{Representation of CW-complexes}\label{SECrep}
A first issue to address in any machine computation of (\ref{EQcohom})
is the machine representation of the CW-complex $X$. Let us briefly consider three possible representations, and then opt for a fourth which seems  better suited to our needs. Consider 
 the torus $\mathbb T=\mathbb S^1\times \mathbb S^1$ with CW-structure involving a
single $0$-cell, two $1$-cells, and a single $2$-cell obtained by taking the product of two copies of the CW-complex $S''$ shown in Figure \ref{fig2} (right). One possibility is to represent
$\mathbb T$ as a free presentation of its fundamental group, involving 2 generators and 1 relator. This representation is readily adapted to $2$-dimensional
CW-complexes with more than one $0$-cell, and is used by 
\cite{MR1743389}
in their algorithm for computing fundamental groups and finite covers of $2$-dimensional cell complexes. However, the  notion of a group presentation does not readily  adapt
 to higher dimensional CW-complexes such as the product $\mathbb S^1\times \mathbb S^1\times \mathbb S^1\times \mathbb S^1$ of four copies of the
CW-complex $S''$.
A second possibility is to subdivide the CW-structure on $\mathbb T = \mathbb S^1\times \mathbb S^1$ in a way that
produces  a simplicial complex. It is possible to triangulate $\mathbb T$ using $7$
vertices, $21$ edges and $14$ triangular faces; the resulting simplicial complex can be stored as a collection of subsets of the vertex set. The approach readily
generalizes to CW-complexes of arbitrary dimension, but the number of
simplices can become prohibitively large. For instance, the $n$-sphere $\mathbb S^n$
admits a CW-structure
with just $2$ cells, whereas any homotopy equivalent simplicial complex
requires at least $2^{n+2}-2$ simplices. A third possibility is to obtain
smaller  simplicial cell structures 
 by relaxing the requirements of a simplicial complex to those
 of a simplicial set. The  $n$-sphere can be represented
as a simplicial set with just two non-degenerate cells.  Simplicial sets,
their fundamental groups, and their (co)homology with trivial coefficients
 have been 
 implemented by John Palmieri in \textsf{SAGE}
 \cite{palmieri}.
We  opt against using simplicial sets as the setting for our algorithms
because it seems to be   
 a non-trivial task to find  small simplicial set representations of some of the 
CW-complexes of interest to us, such as  CW-complexes arising as the complements of knotted surfaces in $\mathbb S^4$. In this paper
 we  opt 
to represent a CW-complex $X$ as a {\em regular} CW-complex $Y$ ({\em i.e.} ones whose attaching maps retrict to homeomorphisms on cell boundaries) together with
a  homotopy equivalence $Y\simeq X$.
 More specifically, we work
with a regular CW-complex $Y$, rather than a simplicial complex,
and  endow it with an
{\em admissible discrete vector field} whose {\em critical cells} are
in one-one correspondence with the cells of $X$.  Definitions of italicized
terms are recalled in Section \ref{SECchn} below. The {\em arrows} in a discrete vector field represent
elementary simple homotopies, and can be viewed as   analogues
of the degeneracy maps of a simplicial set.

\section{Limitations of naive computations}\label{SECnaive}
\begin{figure}
\centerline{\includegraphics[height=4cm]{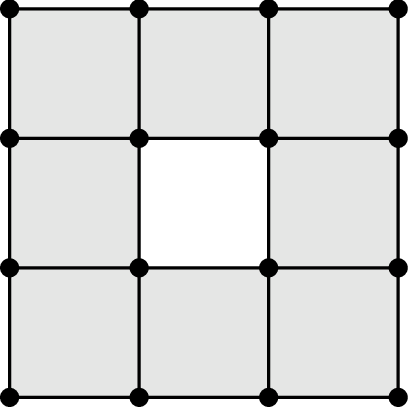}}
\caption{A regular CW-complex $S$ homotopy equivalent to the circle $\mathbb S^1$} \label{fig1}
\end{figure}

\begin{figure}
\centerline{\includegraphics[height=4cm]{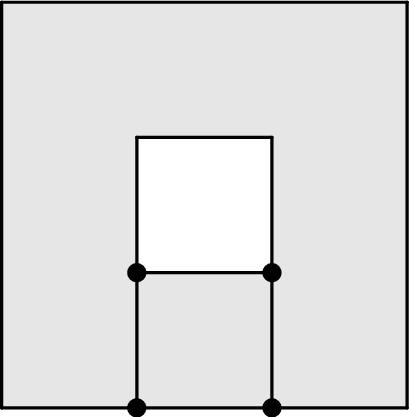}\hspace{2cm}
\includegraphics[height=4cm]{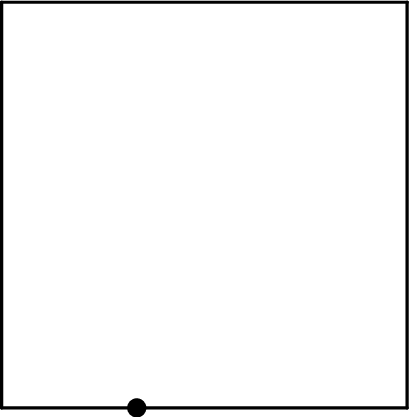}}
\caption{A regular CW-complex $S'$ and a non-regular CW-complex $S''$.  } \label{fig2}
\end{figure}

Having settled on a machine representation for CW-complexes, a second
issue to consider is the size of the computations involved in a naive
implementation of the cellular chain complex of spaces.
Consider the CW-complex $S$ in Figure \ref{fig1} involving 16 vertices, 24
edges, and 8 square $2$-cells.
We may be interested in studying the direct product $Y=S\times S\times S\times S$ which is naturally a regular CW-complex, homotopy equivalent to the $4$-torus, involving a total of $5308416 = 48^4$ cells.
The fundamental group $\pi_1Y$ is free abelian on four generators
$a,b,c,d$ and has a subgroup $H=\langle a^5,b^5,c^5,d\rangle$ of index $125=5^3$. There is a corresponding $125$-fold covering map
$p\colon \widetilde Y_H \rightarrow Y$. We may be interested in the integral homology groups $H_n(\widetilde Y_H,\mathbb Z)$.
 General theory tells us that $\widetilde Y_H$ is homeomorphic to $Y$,
but if we were  ignorant of this general fact then we might consider computing the homology
$H_n(\widetilde Y_H,\mathbb Z) = H_n(C_\ast\widetilde Y_H)$ directly from the
cellular chain complex of $\widetilde Y_H$.
The chain complex $C_\ast\widetilde Y_H$ of the 125-fold cover
can we expressed in terms of the
cellular chain complex
$C_\ast \widetilde Y$ of the universal cover:
\begin{equation}
C_\ast\widetilde Y_H = C_\ast\widetilde Y\,\otimes_{\mathbb ZG} A
\end{equation}
with $G=\pi_1Y$, $A=\mathbb ZG\,\otimes_{\mathbb ZH} \mathbb Z$. The desired
homology
$H_n(\widetilde Y_H,\mathbb Z)$ can be viewed as the homology $H_n(Y,A)$
of $Y$  with local coefficients in the  module $A$.
The chain complex $C_\ast\widetilde Y_H$ is a complex of free abelian groups
of the following form.
\begin{equation}
\begin{minipage}{0.9\textwidth}
$
\rightarrow 0
\rightarrow \mathbb Z^{512000}
\rightarrow \mathbb Z^{6144000}
\rightarrow \mathbb Z^{31744000}
\rightarrow \mathbb Z^{92160000}
\rightarrow \mathbb Z^{164352000}
\rightarrow$

\medskip
\hfill{$\mathbb Z^{184320000}
\rightarrow \mathbb Z^{126976000}
\rightarrow \mathbb Z^{49152000}
\rightarrow \mathbb Z^{8192000}
$}
\end{minipage}\label{EQchncmplx}
\end{equation}
It is not practical to compute $H_n(\widetilde Y_H,\mathbb Z)$ by
applying  the Smith Normal Form algorithm directly to (\ref{EQchncmplx}).
An aim of this paper is to explain how simple homotopy equivalences,
represented as discrete vector fields, can be used to produce
 a smaller  homotopy equivalent chain complex from which the homology can be computed.
An implementation of the techniques is available in the \textsf{HAP} package
\cite{hap} for the \textsf{GAP} system for computational algebra \cite{GAP4},
and is illustrated in the computer code of in Table \ref{tableone}  which
 computes
$H_0(\widetilde Y_H,\mathbb Z)=\mathbb Z$, $H_1(\widetilde Y_H,\mathbb Z)=\mathbb Z^4$,
$H_2(\widetilde Y_H,\mathbb Z)=\mathbb Z^6$,
$H_3(\widetilde Y_H,\mathbb Z)=\mathbb Z^4$,
$H_4(\widetilde Y_H,\mathbb Z)=\mathbb Z$ from the regular CW-complex $Y=S\times S\times S\times S$.
\begin{table}[h]
\begin{verbatim}
gap> A:=[[1,1,1],[1,0,1],[1,1,1]];;
gap> S:=PureCubicalComplex(A);;
gap> S:=RegularCWComplex(S);;
gap> Y:=DirectProduct(S,S,S,S);
Regular CW-complex of dimension 8
gap> Size(Y);
5308416
gap> C:=ChainComplexOfUniversalCover(Y);
Equivariant chain complex of dimension 4
gap> G:=C!.group;;  
[ f1, f2, f3, f4 ]
gap> H:=Group(G.1^5,G.2^5,G.3^5,G.4);
Group([ f1^5, f2^5, f3^5, f4 ])
gap> D:=TensorWithIntegersOverSubgroup(C,H);
Chain complex of length 4 in characteristic 0 .
gap> Homology(D,0);
[ 0 ]
gap> Homology(D,1);
[ 0, 0, 0, 0 ]
gap> Homology(D,2);
[ 0, 0, 0, 0, 0, 0 ]
gap> Homology(D,3);
[ 0, 0, 0, 0 ]
gap> Homology(D,4);
[ 0 ]
\end{verbatim}
\caption{Computing the homology of a 125-fold cover of a regular CW-complex $Y$ with 5308416 cells. }
\label{tableone}
\end{table}
This computation involves a non-regular CW-complex $\widetilde X_H$ which
is homotopy equivalent to
$\widetilde Y_H$. More precisely, $\widetilde X_H$ is a finite cover of
$X=S''\times S''\times S''\times S''$ with $S''$ the non-regular
CW-complex of
Figure \ref{fig2} (right). 
The CW-complex $\widetilde X_H$ involves considerably fewer cells that $\widetilde Y_H$.

For some situations one may wish to work
with a smaller regular
CW-complex
 that is homeomorphic, and not just homotopy equivalent,
 to  $\widetilde Y_H$. In the above example the
number of cells in $\widetilde Y_H$ places its direct construction out of
reach of computers with modest CPU and memory.
 But one could try to simplify the cell structure on  $\widetilde Y_H$ so that 
it
involves fewer cells. The obvious approach is to first
simplify the cell structure on $Y$. In the above example we have
$Y=S\times S\times S\times S$ with $S$ the CW-complex of Figure \ref{fig1}.
We could replace $S$ by the homeomorphic regular CW-complex
 $S'$ of Figure \ref{fig2}
(left) and then work with the space
$Y'=S'\times S'\times S'\times S'$ which is homeomorphic to $Y$.
Using an algorithm for simplifying the CW-structure on a regular CW-complex, the
 computer code in Table \ref{tabletwo} illustrates the construction of
 the covering map
$p\colon \widetilde Y'_H \rightarrow Y'$ which maps each cell of $\widetilde Y'_H$ homeomorphically to a cell of $Y'$.
\begin{table}[h]
\begin{verbatim}
gap> S:=SimplifiedComplex(S);
Regular CW-complex of dimension 2
gap> Size(S);
12
gap> Y:=DirectProduct(S,S,S,S);;
gap> U:=UniversalCover(Y);
Equivariant CW-complex of dimension 8
gap> G:=U!.group;;
gap> H:=Group(G.1^5,G.2^5,G.3^5,G.4);;
gap> p:=EquivariantCWComplexToRegularCWMap(U,H);
Map of regular CW-complexes
gap> UH:=Source(p);
Regular CW-complex of dimension 8
gap> Size(UH);
2592000
\end{verbatim}
\caption{Computing a 125-fold covering map $p\colon \widetilde Y'_H \twoheadrightarrow Y'$ with $Y'$ homeomorphic to a regular CW-complex
with 5308416 cells.}
\label{tabletwo}
\end{table}

The data types,  algorithms, and  mathematical theory underlying
the above example computations are explained in Sections \ref{SECchn} and \ref{SECdata}.

\section{Computing chains on the universal cover}\label{SECchn}

Let $Y$ be a connected regular CW-complex with only finitely many cells. 
One algorithm for constructing a finite presentation for the fundamental group $G=\pi_1Y$ was described in \cite{EllisMrozek},  \cite{ellisbook}, details of which are recalled below. It is well-known that, in general, there is no algorithm for solving the word problem in the finitely presented group $G$. But this unsolvability of the word problem is not relevant to our goals. 
 Let $\widetilde Y$ be the universal cover of $Y$ with canonical CW-structure inherited from $Y$. 
 The cellular chain group $C_n\widetilde Y$ is a free $\mathbb ZG$-module whose free generators  correspond to the $n$-cells of $Y$.
 The module $C_n\widetilde Y$ can be represented on a computer by
specifying the number 
 ${\rm rank}_{\mathbb ZG}(C_n\widetilde Y)$ of free generators
and specifying the free presentation 
for $G$. 
We denote by $e^n_i$ an $n$-cell of $Y$, and we denote by   
 $\tilde e^n_i$ some preferred lift of $e^n_i$ in $\widetilde Y$. For $g\in G$ we 
denote by $g \tilde e^n_i$ the n-cell obtained by applying the action of $g$ to
 $\tilde e^n_i$.  We also use $\tilde e^n_i$ to denote a free generator of the $\mathbb ZG$-module $C_n\widetilde Y$.  
The boundary homomorphism $\partial_n\colon C_n\widetilde Y \rightarrow C_{n-1}\widetilde Y$ is $G$-equivariant and can be represented by specifying $\partial_n(\tilde e^n_i)$  for each  free generator of $C_n\widetilde Y$. We now describe how an
 explicit expression for 
$\partial_n(\tilde e^n_i)$ can be computed using induction on $n$.

We view the $1$-skeleton $Y^1$ as a graph. Let $T$ denote some fixed choice of
 maximal tree in $Y^1$. The tree $T$ has the same vertices as $Y^1$ but, typically, fewer edges.
 Each edge $e$ in $Y^1$ admits two possible orientations;  let us arbitrarily fix one choice of  orientation for each edge.  An oriented
edge $e$ in $Y^1\setminus T$ determines an oriented circuit  in $Y^1$ whose image 
consists of the oriented edge $e$ and some oriented edges from $T$. We denote this  oriented circuit by   
$\omega(e)$.  
 The association $e\mapsto \omega(e)$  induces a function
$\omega \colon \{{\rm edges \ in \ }Y^1\setminus T \}\longrightarrow G$. 
This function extends to a function $\omega \colon \{{\rm edges \ in \ }Y^1 \}\longrightarrow G$ which sends oriented edges in $T$ to the trivial element of the fundamental group.

Let $F$ denote the free group on the symbols $\omega(e)$ where $e$ runs over the edges in $Y^1\setminus T$.  Each $2$-cell $e^2$ in $Y^2$ admits two possible orientations;  let us arbitrarily fix one choice of  orientation for each $2$-cell. The boundary of an oriented $2$-cell then
spells a word  in $F$ that represents the trivial word in the fundamental group $G$.   It is well-known \cite{opacb1122188}, \cite{spanier} that $G$
 admits a
 free presentation
with one generator $\omega(e)$ for each edge $e\subset Y^1\setminus T$, and one relator for each $2$-cell in $Y$.

Let us now consider how to compute an expression for the boundary 
$\partial_n(\tilde e^n_i)$ of each free generator of the $\mathbb ZG$-module
 $C_n(\widetilde Y)$. We first consider the case $n=1$.
The $1$-cell $\tilde e^1_i$ in $\widetilde Y$ maps to the $1$-cell $e^1_i$ in $Y$.
Suppose that $e^0_j, e^0_k$ are the two boundary vertices of $e^1_i$, 
and that the orientation of $e^1_i$ corresponds to a direction which starts at $e^0_j$ and
ends at $e^0_k$. Suppose also that in the cellular chain complex, which involves choices of orientation, the boundary homomorphism satisfies
\begin{equation}
\partial_1(e^1_i) = e^0_k - e^0_j\ .
\end{equation} 
Then we set
\begin{equation}
\partial_1(\tilde e^1_i) = g\tilde e^0_k - \tilde e^0_j\ ,\ \ \ {\rm where\ } g=\omega(e^1_i) .
\end{equation}
Suppose now that we have computed $\partial_1(\tilde e^1_i)$,
$\partial_2(\tilde e^2_i)$, \ldots, $\partial_n(\tilde e^n_i)$ for all free
generators in degrees $\le n$.
To compute  $\partial_{n+1}(\tilde e^{n+1}_i)$
we first compute 
\begin{equation}
\partial_{n+1}e^{n+1}_i = \epsilon_1e^n_{i_1} + \epsilon_2e^n_{i_2} + \cdots + 
\epsilon_me^n_{i_m} \label{EQsix}
\end{equation}
with $\epsilon_i=\pm 1$.
Having computed (\ref{EQsix}) it remains to determine representatives of group elements $g_{i_1},\ldots, g_{i_m}$ in the finitely presented group $G$ such that
the element
\begin{equation}
\partial_{n+1}\tilde e^{n+1}_i = \epsilon_1g_{i_1}\tilde e^n_{i_1} + \epsilon_2g_{i_2}\tilde e^n_{i_2} + \cdots + \epsilon_mg_{i_m}\tilde e^n_{i_m} 
\label{EQseven}
\end{equation}
satisfies $\partial_n(\partial_{n+1}\tilde e^{n+1}_i)=0$.
 We set $g_{i_1}=1 \in G$.
Suppose $\partial(g_{i_1}\tilde e^n_{i_1})$ contains $\pm h \tilde e^{n-1}_j$
 as a summand, some $h\in G$. Then for some $i_j\ne i_1$   the boundary
$\partial(\tilde e^n_{i_j})$ must contain $\mp h' \tilde e^{n-1}_j$
as a summand for some $h'\in G$.
In (\ref{EQseven}) we set $g_{i_j} = hh'^{-1}$. 
 We  continue with this method of matching summands in the boundaries 
$\partial \tilde e^n_{i_k}$ in order to determine all the $g_{i_k}$ in 
(\ref{EQseven}). Note that the method does not require a solution to the word problem in the finitely generated group $G$.

Suppose given  a finitely generated abelian group $A$ and group homomorphism 
$\phi\colon G\rightarrow Aut(A)$ specified on  generators of  $G$. 
This data constitutes a $\mathbb ZG$-module $A$. 
 Starting from our computer representation of the chain complex 
$C_\ast\widetilde Y$, it  is routine to  implement the cochain complex 
$\Hom_{\mathbb ZG}(C_\ast \widetilde Y , A)$ of finitely generated abelian 
groups, and
also the chain complex 
$C_\ast \widetilde Y \otimes_{\mathbb ZG} A$. 
 In particular, given a finite set of elements in $G$ that generates
 a finite index subgroup $H<G$, we can use coset enumeration to construct 
the $\mathbb ZG$-module
$A=\mathbb ZG\otimes_{\mathbb ZH} \mathbb Z$; in this case the chain complex
$C_\ast \widetilde Y \otimes_{\mathbb ZG} A$ can be viewed as 
the chain complex
$C_\ast \widetilde Y_H$ of the finite covering space 
$p\colon \widetilde Y_H \twoheadrightarrow Y$ 
 with $p$ inducing an isomorphism $\pi_1\widetilde Y_H \cong H$. 
 Since $\widetilde Y_H$ is a regular CW-complex, its CW-structure ({\em i.e.}~its face lattice) is completely determined by the chain complex $C_\ast \widetilde Y_H$.
Thus, in principle, we  have an algorithm for computing the finite cover $\widetilde Y_H$.
Furthermore, in principle, the Smith Normal Form algorithm can  be used to compute the local cohomology $H^n(Y,A)$ and homology $H_n(Y,A)$ for any $\mathbb ZG$-module $A$ which is finitely generated over $\mathbb Z$.

However, for practical computations, such as those in Tables \ref{tableone} 
and \ref{tabletwo}, there are two  issues with the above theoretical 
algorithm. 
Firstly, if the above method is applied directly then the
 finitely presented fundamental group $G=\pi_1Y$ will typically involve  
excessive numbers of generators and relators. These numbers need to be 
reduced in a way that retains the relationship between the finite 
presentation and the cellular structure of $Y$ if, for instance,
 one wants to apply coset enumeration and the Reidemeister-Schreier algorithm to  list free presentations of subgroups $H\le G$ of given index $k$ in order to enumerate all $k$-fold covers of $Y$. Table \ref{FIGxt} provides a computation
in which such an enumeration is needed. Secondly, the number of cells
in $Y$ will typically be large, and consequently the number of generators
 of the  chain 
groups 
$C_n\widetilde Y\otimes_{\mathbb ZG} A$  will typically make it impractical 
to apply the Smith Normal Form algorithm directly to the chain complex $C_\ast\widetilde Y\otimes_{\mathbb ZG} A$ or cochain complex 
$\Hom_{\mathbb ZG}(C_\ast \widetilde Y , A)$.
We need a method for reducing the number of cells in $Y$ and $\widetilde Y$ in a way that
 retains the homotopy types of $Y$ and $\widetilde Y$ or even, for some purposes, their homeomorphism types. The code in Table \ref{tableone} reduces the number of cells while retaining only the homotopy type of $Y$. The code in Table \ref{tabletwo} retains the homeomorphism type of $Y$.
 We now describe our approach to addressing these  practical issues.

\begin{figure}[h]
$$\includegraphics[height=4cm]{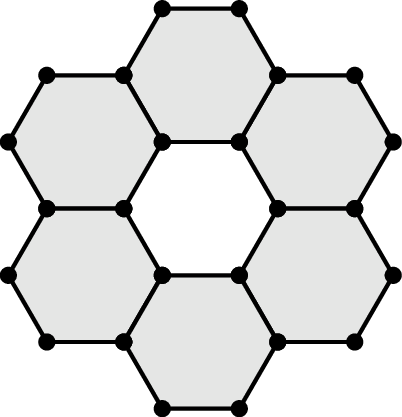}
~~~~~~~~~~~~~~~~~~~~~\includegraphics[height=4cm]{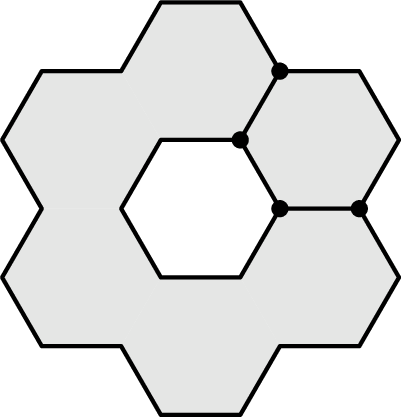}
$$
\caption[Simplification of a regular CW-complex.]{ Regular CW-complex before and after simplification.
}\label{FIGsimplify}
\end{figure}

Let us first recall from \cite{ellisbook} a simple method aimed at reducing the number of cells in a regular CW-complex $Y$ while retaining the homeomorphism type of $Y$.
The simplification procedure invoked in  the computer code of Table 
\ref{tabletwo}
 is based on
the observation that if a regular CW-complex $Y$ contains a $k$-cell $e^k$
 lying
in the boundary of precisely two $(k+1)$-cells $e^{k+1}_1, e^{k+1}_2$ with identical coboundaries then
these three cells can be removed and replaced by a single cell of dimension
 $k+1$. The topological space $Y$ is unchanged; only its CW-structure changes.
 The resulting CW-structure will not in general be regular. However,
it will be regular if the sets  $V_0, V_1, V_2$ of cells lying
in the boundaries of $e^k, e^{k+1}_1, e^{k+1}_2$ respectively are such that
$V_1\cap V_2=V_0\cup\{e^k\}$. The result of this kind of
 simplification procedure is illustrated in Figure \ref{FIGsimplify}.
 The procedure is formalized as Algorithm \ref{ALGsimplify}.
\begin{algorithm}
\caption{Simplification of a regular CW-structure.}
\label{ALGsimplify}
\begin{algorithmic}[1]
\Require A regular CW-complex $Y$.
\Ensure A regular CW-complex $X$ with $Size(X)\le Size(Y)$ and with $|X|$ homeomorphic to $|Y|$.
\Procedure{}{}
\State Let $X$ be a copy of $Y$.

\While{there exists a cell
 $e^k$ in $X$ with precisely two cells $e^{k+1}_1, e^{k+1}_2$ in its
coboundary and these two cells have  identical coboundaries}
\State Compute the sets $V_0$, $V_1$, $V_2$ of cells in
the boundaries of $e^k$, $e^{k+1}_1$, $e^{k+1}_2$.
\If{$|V_1\cap V_2|=1+|V_0|$}

\State Remove the cells $e^k$,  $e^{k+1}_1$, $e^{k+1}_2$ from $X$
and add a new
$(k+1)$-cell $f^{k+1}$ to $X$
 whose boundary is the union of the boundaries of $e^{k+1}_1$ and $e^{k+1}_2$ minus the cell $e^k$. Adjust coboundaries accordingly.
\EndIf
\EndWhile
\State Return $X$.
\EndProcedure
\end{algorithmic}
\end{algorithm}
It can be useful, for instance, in the study of the homeomorphism type
 of $CW$-complexes $Y$ arising as knot complements in $\mathbb R^3$, especially when $Y$ is constructed from experimental data, such as knotted proteins, with 
  natural CW-structure involving a vast number of cells (see \cite{EllisMrozek}).

We now turn to a method for reducing the number of cells in $Y$  in a way that
 retains only the homotopy type of $Y$. For this we first recall the following
notion.

\begin{definition}\label{DEFdiscretevectorfield}
A {\em discrete vector field} on a regular CW-complex $Y$ is a collection of pairs $(s , t)$, which we call {\em arrows} and denote by $s \rightarrow t$,  
 satisfying 
\begin{enumerate}
\item  $s,t$ are cells of $X$ with $\dim(t) =\dim(s) +1$ and with
$s$ lying in   the boundary of $t$. We say that $s$ and $t$ are {\it involved}
 in the arrow, that $s$ is the {\it source} of the arrow, and that
$t$ is the {\it target} of the arrow.

\item any cell is involved in at most one arrow.

\end{enumerate}

\end{definition}
The term {\em discrete vector field} is due to  \cite{MR1612391}.
 In an earlier work \cite{MR968920} Jones calls this concept  a
 {\it marking}.

A discrete vector field is
said to be {\em finite} if it consists of just finitely many arrows.
By a {\it chain}\index{chain} in a discrete vector field we mean a sequence of
arrows 
$$\ldots , s_1 \rightarrow t_1, s_2 \rightarrow t_2, s_3
\rightarrow t_3, \ldots $$  where $s_{i+1}$
lies in the boundary of $t_i$ for each $i$.
A chain is  a {\em circuit} \index{circuit} if it is of finite length with source $s_1$
of the  initial arrow $s_1\rightarrow t_1$ lying in the boundary of  the
target $t_n$ of the final arrow $s_n\rightarrow t_n$.
A discrete vector field is
said to be {\it admissible} \index{admissible} if it contains no circuits and
  no chains that extend
infinitely to the right. When the  CW-complex $Y$ is finite 
  a discrete
vector field is admissible if it contains no circuits. We say that an
admissible discrete vector field is {\it maximal}\index{maximal discrete vector field} if it is not possible to
add an arrow while retaining admissibility. A cell in $X$ is said to be
{\it critical} \index{critical cell} if it is not involved in any arrow.
\begin{theorem}\cite{MR1939695,MR1612391}\label{THMbasic}
 If $Y$ is a regular CW-complex with admissible discrete vector field then there is a homotopy equivalence
$$Y\simeq X$$
where $X$ is a CW-complex whose cells are in one-one correspondence with
  the critical cells  of $X$.
\end{theorem}
An arrow on $Y$ can be viewed as representing a simple homotopy, as introduced in \cite{MR0035437}. The theorem just says that an admissible discrete vector field represents some sequence of (in some sense 'composable') simple homotopies  starting at $Y$ and ending at $X$. 

There are various algorithms for constructing an admissible discrete vector field on a finite regular CW-complex $Y$. The algorithm implemented in the {\sf HAP} package  is  Algorithm \ref{ALGdiscretevectorfield}, which we recall from \cite{ELLISHEGARTY}, \cite{ellisbook}. 
\begin{algorithm}[h]
\caption{Discrete vector field on a regular CW-complex.}
\label{ALGdiscretevectorfield}

\begin{algorithmic}[1]

\Require A finite regular CW-complex $Y$

\Ensure A  maximal admissible discrete vector field on $Y$.

\Procedure{}{}

\State Partially order the cells of $Y$ in any fashion.

\State At any stage of the algorithm each cell  will have precisely one of the following three states: (i) {\em critical}, (ii)
{\em potentially critical}, (iii) {\em non-critical}.

\State Initially deem all cells of $Y$ to be potentially critical.

\While{there exists a potentially critical cell}

\While{there exists a pair of potentially critical cells $s,t$
  such
 that:  $\dim(t)=\dim(s)+1$;  $s$ lies in the boundary of $t$;  no other
potentially critical cell of dimension $\dim(s)$ lies in the boundary of $t$;}

\State Choose such a pair $(s,t)$ with $s$ minimal in the given partial ordering.

\State Add the arrow $s\rightarrow t$ and deem $s$ and $t$ to be  non-critical.

\EndWhile

\If{there exists a potentially critical cell}

\State Choose a minimal potentially critical cell and deem it
 to be  critical.

\EndIf

\EndWhile

\EndProcedure
\end{algorithmic}
\end{algorithm}
In line 2 of
the algorithm the cells could be partially ordered in some way
 that ensures any cell of  dimension $k$ is less than all cells of
 dimension $k+1$. This partial ordering  guarantees that the resulting discrete
 vector field on a path-connected regular CW-complex $Y$ will have a unique
critical $0$-cell. 

A CW-complex $X$ is said to be {\em reduced} if it has only one cell in dimension $0$. There is a standard correspondence between the $2$-skeleton of such a space and a presentation of its fundamental group $G=\pi_1X$. This correspondence 
together with Theorem \ref{THMbasic} and Algorithm 
\ref{ALGdiscretevectorfield} constitute our algorithm for finding a 
presentation for the fundamental group of any connected CW-complex $Y$. The 
presentation has one generator for each critical $1$-cell of $Y$ and one 
relator for each critical $2$-cell.  
Each oriented critical $2$-cell $e \subset Y$ determines an oriented circuit
$\omega(e)$ in $Y^1$; the discrete vector field provides a deformation of this circuit $\omega(e)$ into a circuit $\omega'(e)$ each of whose oriented
edges $t$
is either critical or else the target of some arrow $s\rightarrow t$ with $s$ a $0$-cell. 
The subsequence of $\omega'(e)$ consisting of the oriented critical edges spell a word in the free group on the generators of the presentation. For illustrations and further details of this algorithm see \cite{EllisMrozek}, \cite{ellisbook}.

Let us now address the question of how to reduce the number of cells in  
the universal cover $\widetilde Y$ of a connected regular CW-complex in a way 
that
 retains its homotopy type.  It suffices to note 
that any discrete vector field on $Y$ induces a discrete vector field on 
$\widetilde Y$: there is an arrow $\tilde s \rightarrow \tilde t$ on
$\widetilde Y$ if and only if if $s\rightarrow t$ is an arrow on $Y$ where
  $\tilde s, \tilde t$ are cells in the universal cover that map to $s,t$ , and  where $\tilde s$ lies 
in the boundary of $\tilde t$. Furthermore, if the discrete vector field on $Y$ is admissible then so too is the induced vector field.  Let $X\simeq Y$ be the homotopy equivalence of Theorem \ref{THMbasic}. This equivalence induces a homotopy equivalence of universal covers $\widetilde X \simeq \widetilde Y$
and a chain homotopy equivalence of cellular chain complexes
$C_\ast\widetilde X \simeq C_\ast\widetilde Y$. The cells in $\widetilde X$  are in one-one correspondence with the critical cells in
the induced discrete vector field on $\widetilde Y$. The cellular chain complex $C_\ast\widetilde X$ on the (typically non-regular) CW-complex $X$ is readily constructed directly from the cell structure and admissible discrete vector field of $Y$. This construction is implemented in {\sf HAP} as the function
{\tt ChainComplexOfUniversalCover(Y)} which inputs $Y$, constructs a maximal discrete vector field on $Y$, uses this vector field to compute a presentation
for $\pi_1Y \cong \pi_1X$, and finally returns the $\pi_1Y$-equivariant chain complex
$C_\ast(\widetilde X)$. 

\section{From diagrams to CW-complexes}\label{SECprojection}
Suppose that $N$ is a compact topological submanifold  of $\mathbb R^{n+2}$ with $n=1$ or $2$, 
and that $p\colon \mathbb R^{n+2} \twoheadrightarrow \mathbb R^{n+1}$ is a 
projection onto a hyperplane. Let $M \subset \mathbb R^{n+2}$ denote a 
submanifold homeomorphic to the closed unit disk, whose interior $\mathring M$ 
contains $N$.    
Set $D^{n+1}=p(M)$.  We say that the  pair $p(N) \subset D^{n+1}$ is a {\em diagram} for $N$. Under fairly general conditions it is possible to embellish the diagram with  information on the preimage $p^{-1}(y)$ of certain
 points  $y\in N$ 
 so that the ambient isotopy type of the embedding $N\hookrightarrow \mathbb R^{n+1}$ can be recovered from the embellishment. We are then interested in using 
the embellished diagram to construct a regular CW-structure on $M$  containing
 a CW-subspace  ambient isotopic to  $N$.   
 We illustrate the idea with several examples.

\subsection{A closed $1$-manifold $N$ in $\mathbb R^3$.}
\begin{figure}[h]
\centerline{\includegraphics[height=3cm]{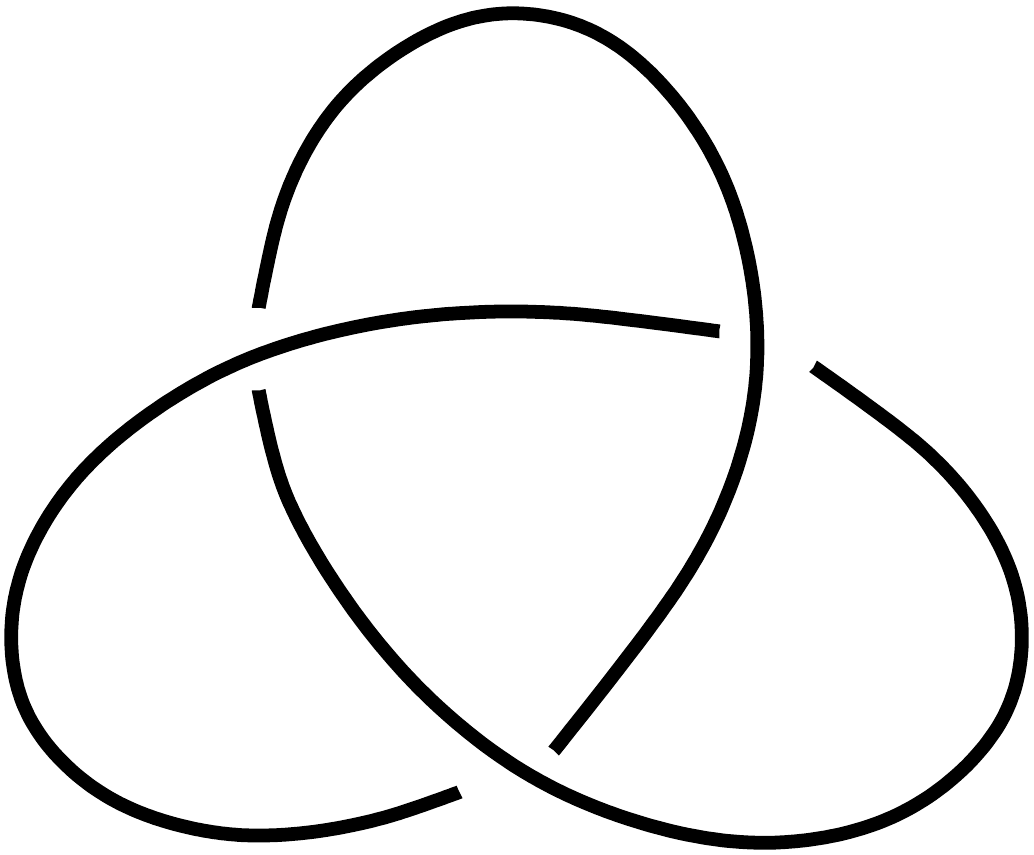}~~~~~ 
\includegraphics[height=3cm]{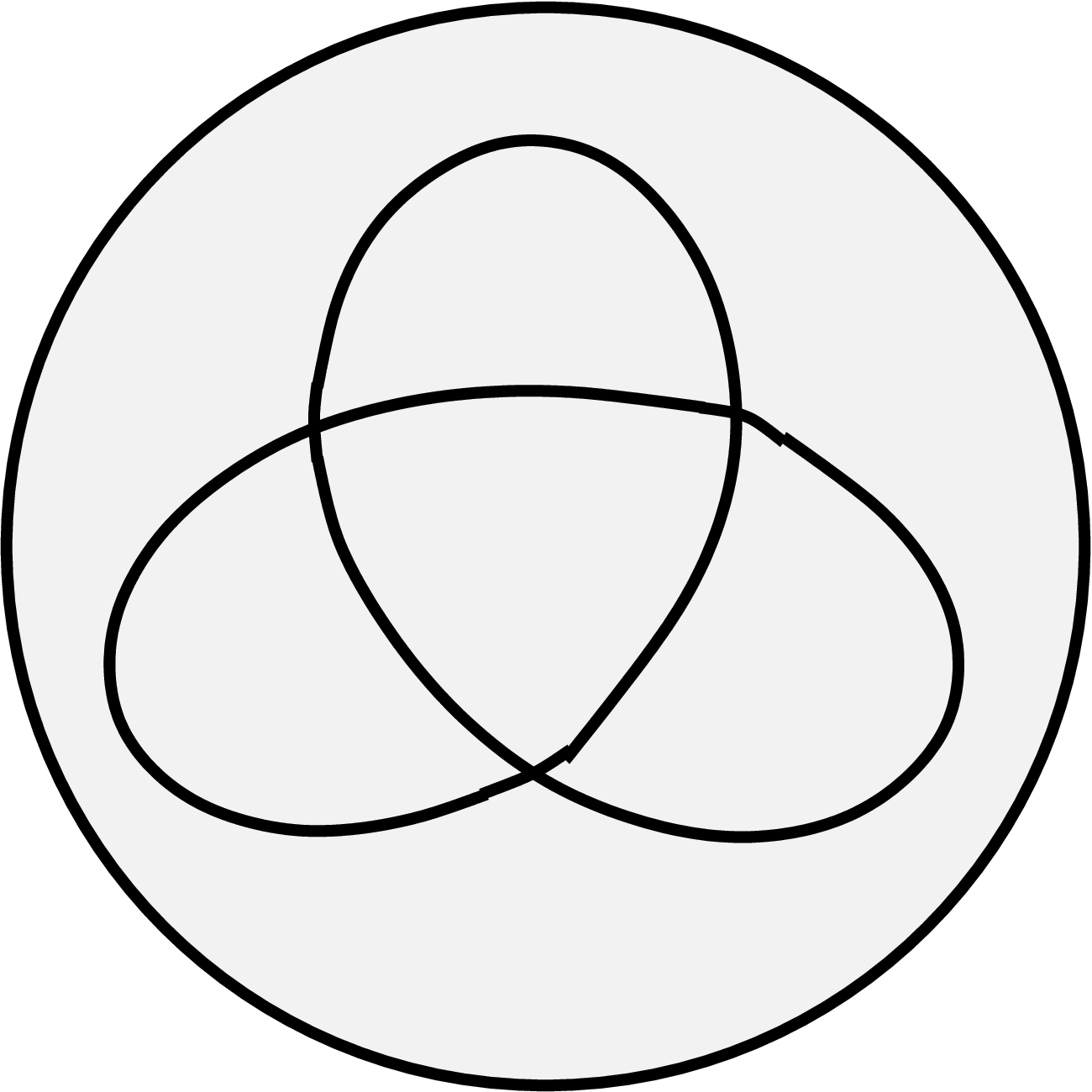} ~~~~~ 
\includegraphics[height=3cm]{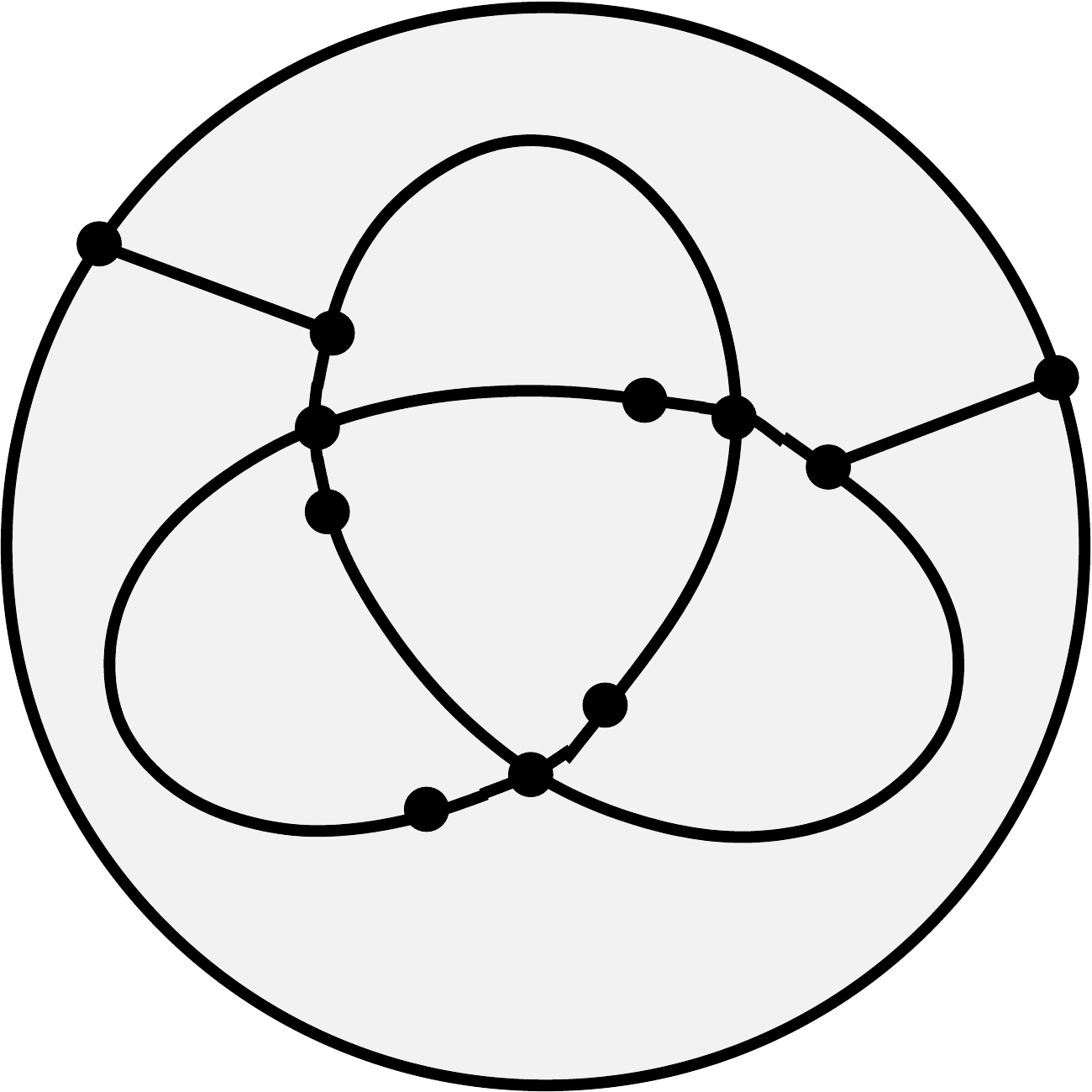}}
\caption{Trefoil knot.}\label{FIGtreff}
\end{figure}

As an illustration let us consider the trefoil knot whose embedding 
$\kappa\colon \mathbb S^1 \hookrightarrow \mathbb R^3$ is illustrated in 
Figure \ref{FIGtreff} (left). We set $N$ equal to the image of $\kappa$.  
A diagram $(D^2,p( N))$ for the trefoil is shown in Figure 
\ref{FIGtreff} (middle). A CW-structure on $D^2$ is shown in Figure 
\ref{FIGtreff} (right). The CW-structure is designed to reflect the under crossing/over crossing structure of the knot projection. Let the unit interval
$[0,1]$ be given a CW-structure  involving two $0$-cells and one $1$-cell. The 
\begin{figure}[h]
\centerline{\includegraphics[height=3cm]{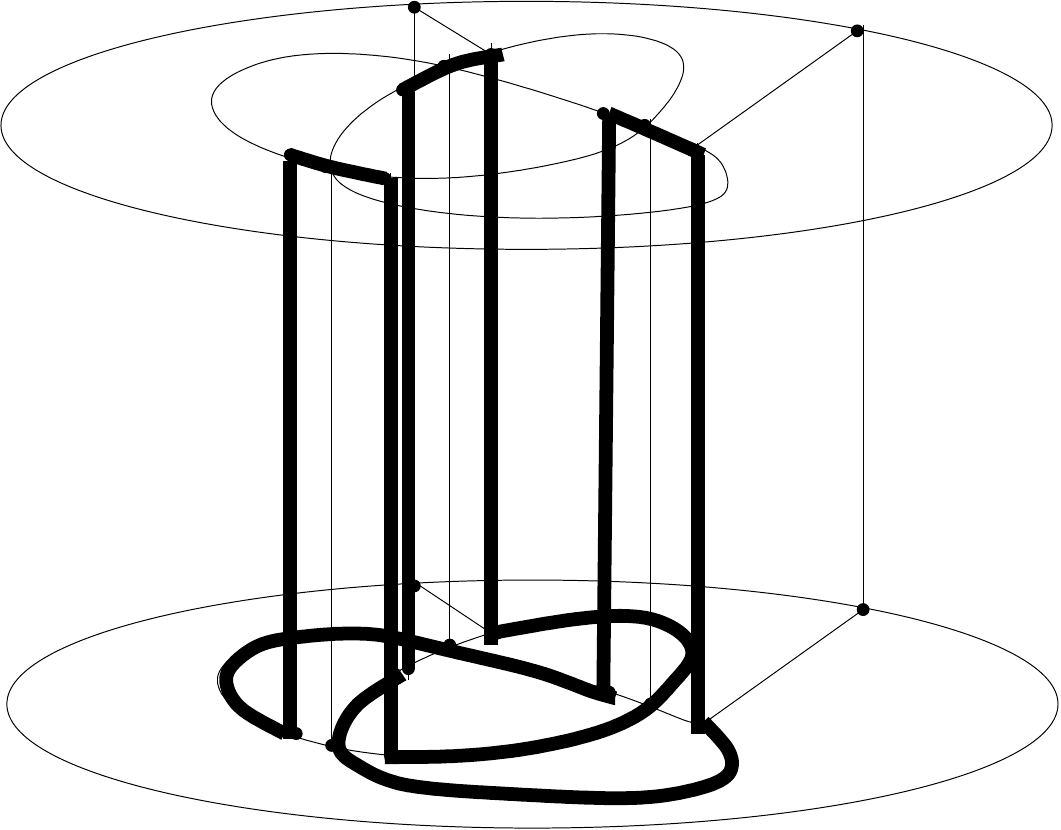}}
\caption{CW-subcomplex  of $M$ ambient isotopic to the Trefoil knot.}\label{FIGtreffknot}
\end{figure}
CW-complex formed by adjoining a $2$-cell and $3$-cell to each end  of
 the cylinder $D^2\times [0,1]$  is homeomorphic to a $3$-ball $M$ and contains, in its interior, a CW-subcomplex $N$ ambient isotopic to the trefoil knot (see Figure \ref{FIGtreffknot}). 
The construction of this CW-structure on the trefoil embedding readily
extends to an algorithm which inputs a symbolic representation of a knot or link and returns a CW-decomposition of a $3$-ball $M$ containing the knot or link as a subcomplex of the $1$-skeleton. 

We are interested in studying the complement $M \setminus N$. However, the difference of CW-complexes is not a CW-complex. This technicality can be overcome by constructing a small open tubular neighbourhood $N_\epsilon$ of $N$ and extending the cell structure on $M\setminus N$ to a CW-structure on $M\setminus N_\epsilon$. Details of the construction are given in Section \ref{SECnbhd}.

\subsection{A $3$-manifold  $N$ in $\mathbb R^3$, and the complement of its interior.}\label{SUB3man}
\begin{figure}
\centerline{\includegraphics[height=5cm]{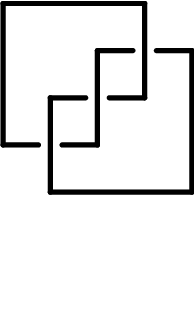}
\hspace{1cm}\includegraphics[height=6cm]{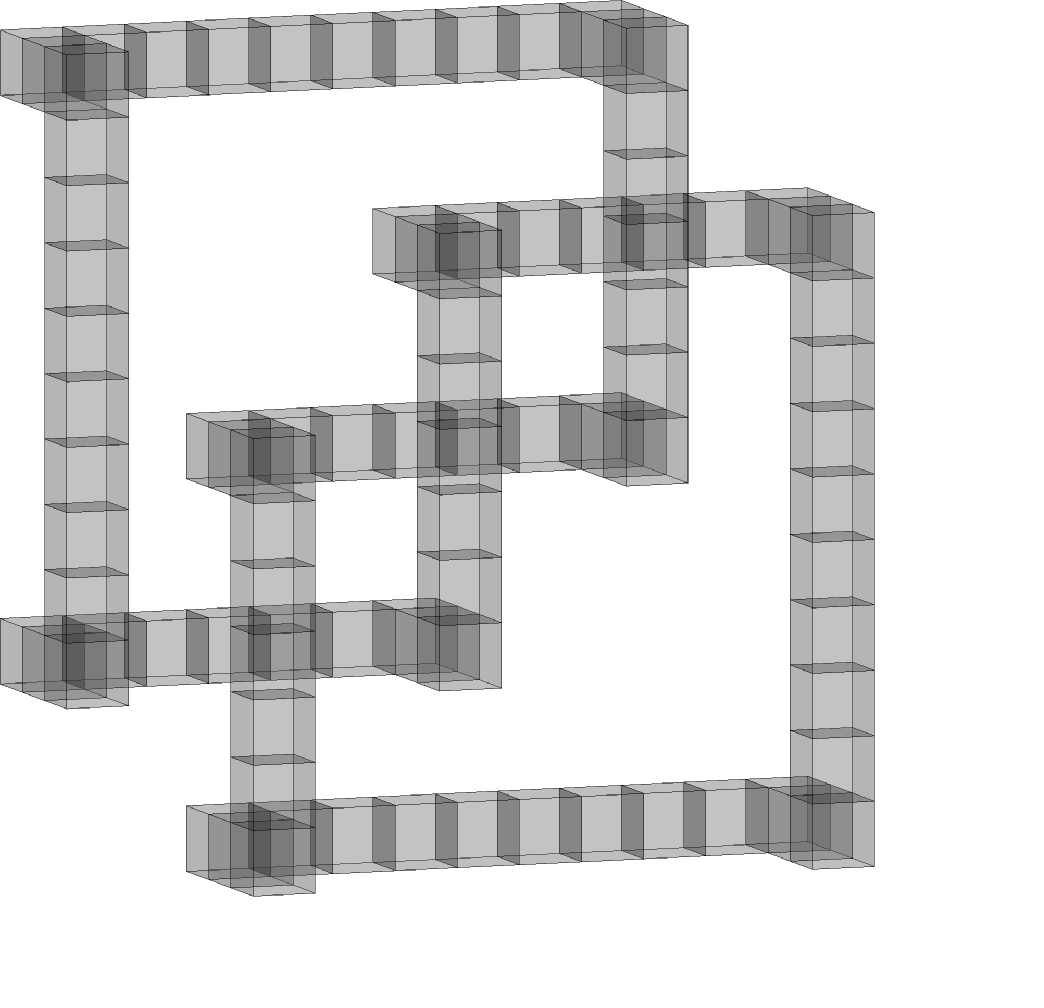}
}
\caption{Arc diagram for the trefoil knot (left) and corresponding pure cubical complex (right).}\label{FIGarc}
\end{figure}

There are a number of succinct symbolic representations of links in $\mathbb R^3$. One representation is based on  arc diagrams of links such as the arc diagram for the trefoil shown in Figure \ref{FIGarc} (left).
An {\em arc diagram} is a planar link diagram involving only horizontal and vertical straight lines, with vertical lines always passing over the horizontal lines.
 An arc diagram can be represented symbolically by a suitable
 ordered list $\ell$
of ordered pairs of positive integers. The list $\ell=[[2,5],[1,3],[2,4],[3,5],[1,4]]$
represents the arc diagram of Figure \ref{FIGarc}.
The $i$th pair in the list  specifies the starting column and end
 column of the $i$th horizontal line, where the bottom line is the $1$st 
horizontal line, and the leftmost column is the $1$st column.

An arc diagram  leads to a very natural, though somewhat inefficient,
  representation of the corresponding
link as a $3$-dimensional CW-manifold $N$ with boundary. 
To explain how, we  endow $\mathbb R^n$ with the regular
 CW-structure in which each $n$-cell has closure equal to a unit cube with standard CW-structure and with centre  an integer vector $c\in \mathbb Z^n$. A finite CW-subcomplex of $\mathbb R^n$ is said to be an $n$-dimensional {\em cubical complex}. Such a complex is said to be {\em pure} if each cell lies in the closure of some $n$-cell.
An arc presentation for $N$  gives rise in a fairly obvious fashion
 to a diagram $(D^2,p(N))$ with $D^2$ a 
$2$-dimensional pure cubical complex homeomorphic to a disk. Now let the 
real interval $[0,5]$ be given a 
CW-structure with integers the $0$-cells.  The desired 
$3$-dimensional pure cubical complex $N$ is realized as a CW-subcomplex of the 
direct product $D^2\times [0,5]$.

The $3$-dimensional pure cubical complex $N_\ell$  corresponding to the  example
list
 for the trefoil is shown in Figure \ref{FIGarc} (right). Let $M$ denote a contractible pure cubical complex in $\mathbb R^3$ whose interior contains $N_\ell$.
Then
the complement $X=M\setminus \mathring{N_\ell}$ is a pure cubical complex homotopy equivalent to the complement  $\mathbb R^3\setminus N_\ell$.
 The construction of the CW-complex $X$ is implemented in {\sf HAP}, with $M$
chosen to be a minimal solid rectangular pure cubical complex whose interior contains $N_\ell$.
For the trefoil knot of Figure \ref{FIGarc} (right) the implemented CW-complex $X$ contains 13291 cells.

\subsection{A more efficient complement of the interior of a $3$-manifold $N$ in $\mathbb R^3$.} \label{SUBefficient}
Pure cubical complexes are particularly
 useful as a tool for converting a range of
 experimental data  into regular CW-complexes in order to analyse underlying
  topological features (see for instance \cite{ellisbook}). However, the resulting
CW-complexes tend to be large and cumbersome
to work with on a computer. From the viewpoint of theoretical knot theory, one
would like  algorithms  for converting a symbolic description of a link
or knotted surface complement into a regular CW-complex with relatively few cells so that computations of fundamental groups and local cohomology run efficiently. 
  We now explain, in four steps, how to construct a smaller CW-structure on the complement  $Y_\ell=M\setminus \mathring{N_\ell}$ of a link $N_\ell$
arising from an
arc presentation $\ell$.

\begin{figure}[h]
\centerline{\includegraphics[height=5cm]{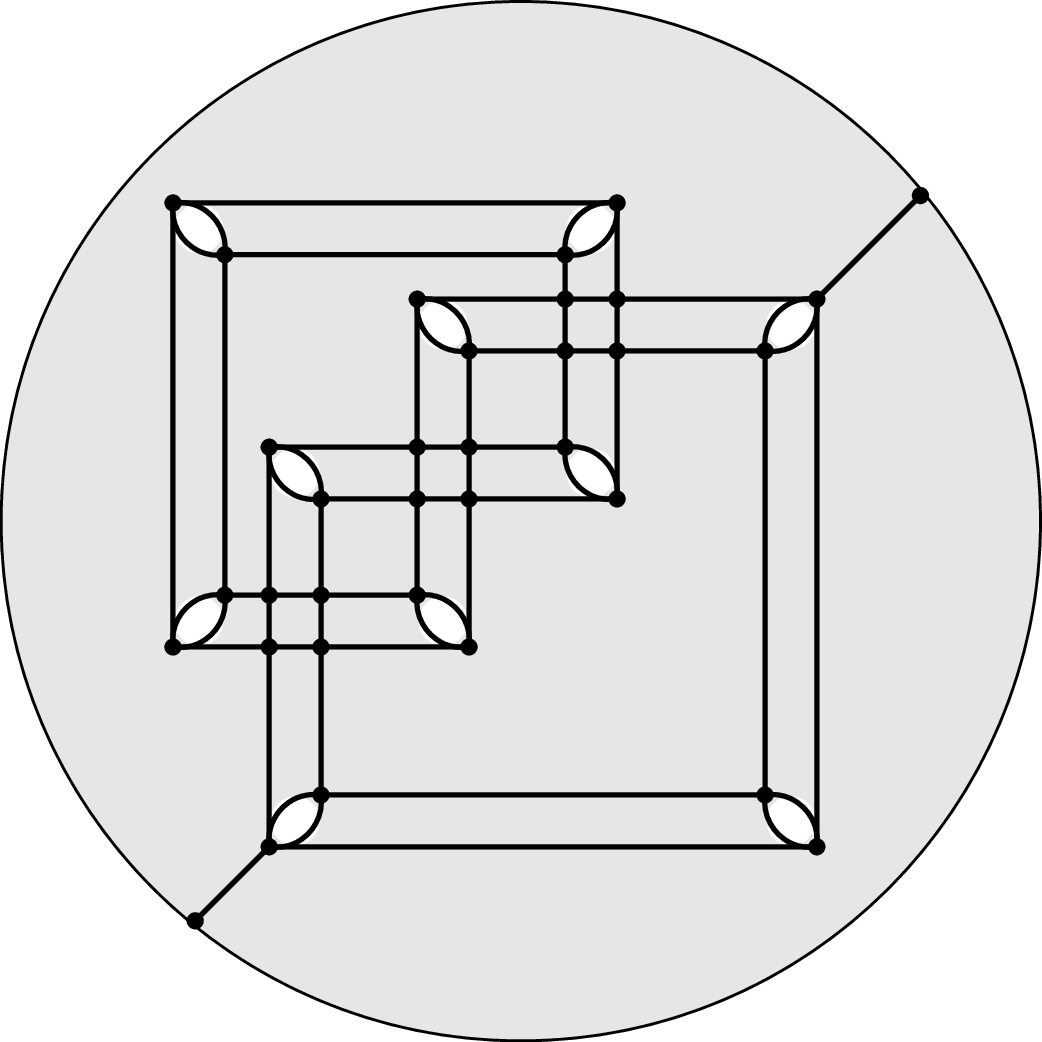} }
\caption{CW-structure on a $2$-dimensional disk $D^2_{\ell}$ with ten holes.}\label{FIGcw}
\end{figure}

Step 1.
Let us suppose that the arc diagram for the list $\ell$ involves precisely
$h$ horizontal
 lines and precisely $k$ crossings. Let $D^2_{\ell}$ denote the unit
$2$-disk with $2h$ non-overlapping subdisks removed.
The arc diagram for $\ell$ then determines a (nearly) canonical CW structure on
$D^2_{\ell}$. 
 This is illustrated in Figure \ref{FIGcw} for the trefoil. (The CW structure could be made canonical by, for instance, insisting that the two $0$-cells on the boundary of the disk be connected to the bottom left-most and top right-most vertices of the inner diagram.)
In general the
 CW structure has $V=4h+4k+2$ vertices, $E=8h+ 8k +4$ edges, and $E-V-2h+1$ cells of dimension $2$. The formula for the number of $2$-cells is derived from the Euler characteristic $\chi(D^2_{\ell}) = 1-2h$.
 The space $D^2_{\ell}$ has a total of $14h  + 16k  +9$ cells.

Step 2. Let $\mathbb I$ denote the unit interval with $CW$-structure involving two vertices and one edge. Now form the direct product $D^2_{\ell}\times \mathbb I$. This is a CW-complex which we view as a solid cylinder from which $2h$ vertical tubes, running from bottom to top, have been removed.
 It involves a total of $3(14h  + 16k  +9)$ cells.

Step 3.
 Corresponding to each horizontal line in the arc diagram glue
 one $2$-cell to the bottom of  $D^2_{\ell}\times \mathbb I$, and corresponding to each vertical line in the arc diagram glue one $2$-cell to the top of
$D^2_{\ell}\times \mathbb I$. More precisely, these $2h$ $2$-cells are glued so that the resulting CW-complex
$W_\ell=D^2_{\ell}\times \mathbb I \cup \bigcup_{i-1}^{2h}e^2$
is homotopy equivalent to the link complement $\mathbb R^3\setminus M_\ell$.

Step 4. Glue one $3$-cell and one $2$-cell to the bottom of $W_\ell$, and
glue one $3$-cell and one $2$-cell to the top of $W_\ell$  in a way that the resulting CW-complex $Y_\ell=W_\ell \cup \bigcup_{i=1}^2 e^2 \cup \bigcup_{i=1}^2 e^3$
is homeomorphic to $M\setminus \mathring N_\ell$.  The CW-complex $Y_\ell$ involves a total of $3(14h  + 16k +9) +2h +4$ cells.

\subsection{A closed $2$-manifold $N$ in $\mathbb R^4$.}
\begin{figure}
$$\includegraphics[height=4cm,angle=90,origin=c]{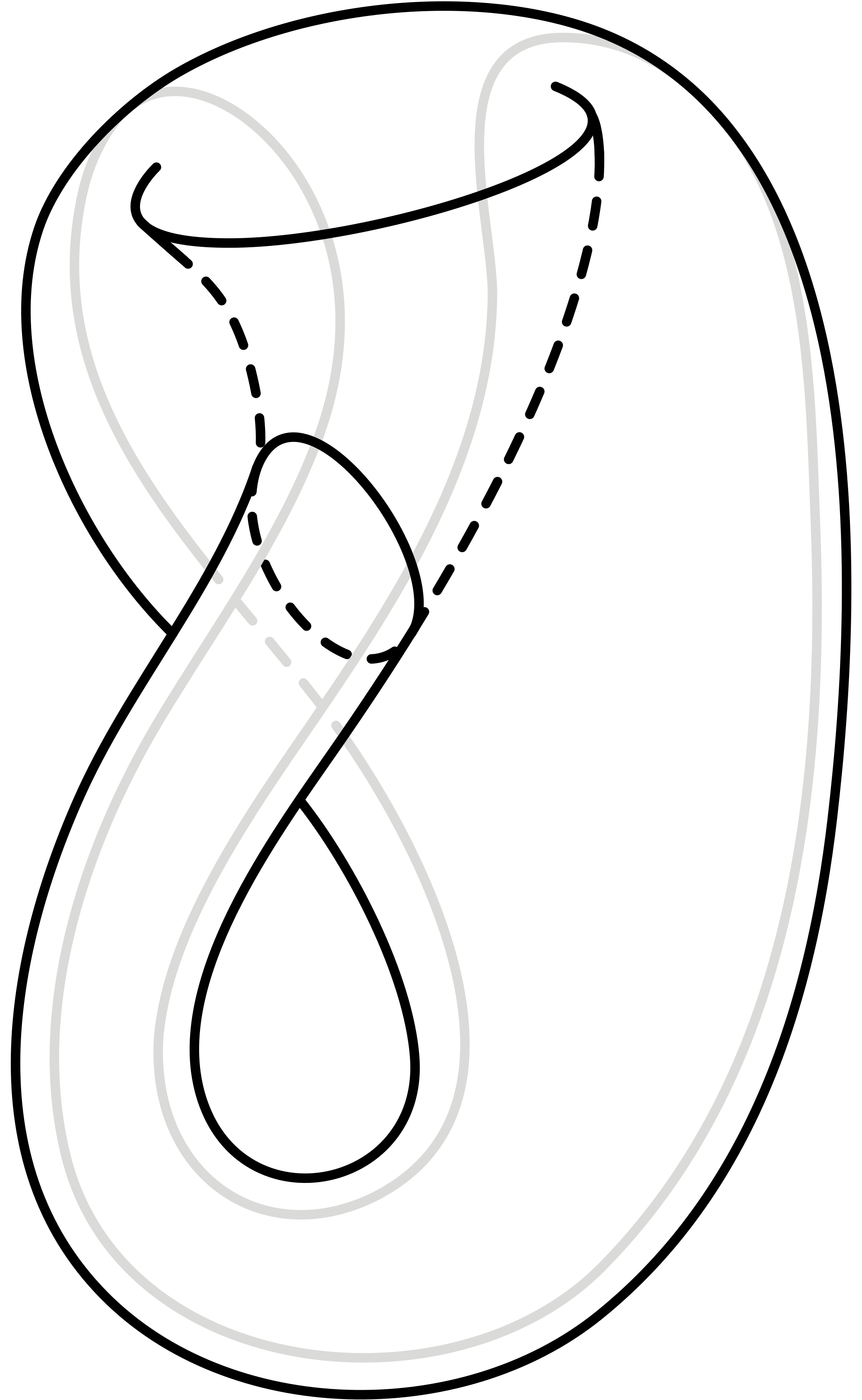}$$
\caption{Klein bottle}\label{FIGklein}
\end{figure}

In order to specify an embedding $N\hookrightarrow \mathbb R^4$  of a closed
surface $N$ we can start by specifying  a diagram $(D^3,p(N))$ with $D^3$ a
 regular 
CW-complex homeomorphic to a closed $3$-ball, and with  $p(N) \subset D^3$ a 
CW-subcomplex 
 lying in the interior of $D^3$. The CW-subcomplex $p(N)$ need not be a surface; we   allow it to have {\em singular points} $y\in p(N)$ for which there
exists  
no open set of $p(N)$  containing  $y$  and  homeomorphic to an open $2$-disk. 
But we do require the collection of singular points to be either empty or to
form a 
closed $1$-manifold.  
For instance, in the usual projection of the Klein bottle into $\mathbb R^3$ the singular points form a circle (see Figure \ref{FIGklein}).
We say that the CW-complex $p(N)$ is a {\em self-intersecting surface}. 
To complete the specification of the embedding of the surface $N$ we 
set $M=D^3\times [0,k]$ where $k$ is a positive integer and the interval is given a
CW-structure with $0$-cells the integers. Let $p\colon M\twoheadrightarrow D^3$
be the projection.
 To specify a CW-manifold $N\subset M$
we just need to specify, in some fashion, the cells in the preimage $p^{-1}(e)$
for each cell $e\subset p(N)$. 

The discussion in Subsection \ref{SUBefficient} can be interpreted as a method
for representing an embedding of one or more tori into $\mathbb R^3$ by means of an arc presentation. This representation is readily
 extended to one for self-intersecting tori whose singular points form a $1$-manifold.

Given an embedding of a closed surface $N$ into $\mathbb R^4$, 
we can then construct a small open tubular neighbourhood $N_\epsilon$ of $N$ and extend the cell structure on $M\setminus N$ to a CW-structure on $M\setminus N_\epsilon$. Details
of the neighbourhood construction   are given in Section \ref{SECnbhd}.

\section{Some data types}\label{SECdata}

For completeness, we recall from \cite{ellisbook} the main data types used in implementing the above discussion.

\begin{datatype} \label{DATAregularcw} \index{data type:: regular CW-complex}
A regular CW-complex $X$ is represented as a component object \textsf{X}
with the following components:
\begin{itemize}
    \item {\sf X!.boundaries[n+1][k]} is a list of integers $[t,a_1,...,a_t]$
 recording that the $a_i$th cell of dimension  $n-1$ lies
in the boundary of the $k$th cell of dimension $n$.
    \item {\sf X!.coboundaries[n][k]} is a list of integers $[t,a_1,...,a_t]$
 recording that the $k$th cell of dimension  $n$ lies in the boundary of the $a_i$th cell of dimension $n+1$.
    \item {\sf X!.nrCells(n)} is a function returning the number of cells in dimension $n$.
    \item {\sf X!.properties} is a list of properties of the complex, each property stored as a pair such as $[``dimension",4]$.
\end{itemize}
\end{datatype}

\begin{datatype}
A discrete vector field on a regular CW-complex is represented as
 a regular CW-complex \textsf{X} with the following additional components, each of which is a $2$-dimensional array:
\begin{itemize}
\item \textsf{X!.vectorField[n][k]} is equal to the integer $k'$ if there is
an arrow from the $k$th cell of dimension $n-1$ to the $k'$th cell of
dimension $n$. Otherwise \textsf{X!.vectorField[n][k]} is unbound.
\item \textsf{X!.inverseVectorField[n][k']} is equal to the integer $k$
 if there is an
arrow  from the $k$th cell of dimension $n-1$ to the $k'$th cell of
dimension $n$. Otherwise \textsf{X!.inverseVectorField[n][k']} is unbound.
\end{itemize}
\end{datatype}

\begin{datatype} \label{DATgcwcomplex}
A    $G$-CW-complex $Y$ with regular CW-structure 
 is represented by a
 component object \textsf{Y}  
 which consists of the
 following components:
\begin{itemize}
\item \textsf{Y!.group} is a group $G$. (This can be stored in one of several
possible formats including: finitely presented group; matrix group;
power-commutator presented group; finite permutation group ... .)
\item \textsf{Y!.elts} is a list of some of the elements of $G$.
\item \textsf{Y!.dimension(n)} is a function that returns a non-negative integer
equal to the number of distinct $G$-orbits of $n$-dimensional cells in the
 $G$-equivariant CW-complex $|Y|$
modelled by the data type. We let $e^n_k$ denote a
fixed representative of the $k$th orbit of $n$-dimensional cells.
\item \textsf{Y!.stabilizer(n,k)} returns the  subgroup of $G$ consisting of those
elements $\phi \in G$ that  map  $e^n_k$ homeomorphically to itself.
\item \textsf{Y!.boundary(n,k)} returns a list $[[i_1,g_1], \ldots, [i_m,g_m]]$ where each term $[i,g]$ is a pair of non-zero integers with $g$ positive.
A pair $[i,g]$ records that $\phi_g(e^{n-1}_{|i|})$ lies in the boundary
of $e^n_k$, where $\phi_g$ is the $g$th term of the list \textsf{Y!.elts}. (The sign of $i$ plays a role  when we consider homology.)

\item \textsf{Y!.properties} is a list of properties of the space, each property
 stored as a pair such as $[``dimension",4]$.
\end{itemize}
\end{datatype}

\begin{datatype} \label{DATAchaincomplex}
A chain complex $C_\ast$ of finitely generated free modules over  $\mathbb Z$
 is represented as a component object \textsf{C}
with the following components:
\begin{itemize}

\item \textsf{C!.dimension(k)} is a function which returns the rank of the module $C_k$.

\item \textsf{C!.boundary(k,j)} is a function which returns the image in $C_{k-1}$
 of the $j$th free generator of $C_k$. Elements in $C_{k-1}$ are represented
 as vectors ({\it i.e.~}lists) over $\mathbb Z$
of length equal to the rank of $C_{k-1}$.

\item \textsf{C!.properties} is a list of properties of the complex, each property stored as a pair such as $[``length",4]$.

\end{itemize}

\end{datatype}

\section{The spun Hopf link and the tube of the welded Hopf link}\label{SEChopf}
We now give a computer proof of a variant on Theorem 10 in \cite{MR2441256}. We opt to work in the inefficient category of pure cubical complexes as a means of demonstrating that computations are practical even when CW-complexes involve large numbers of cells.

\begin{figure}
\centerline{\includegraphics[height=7cm]{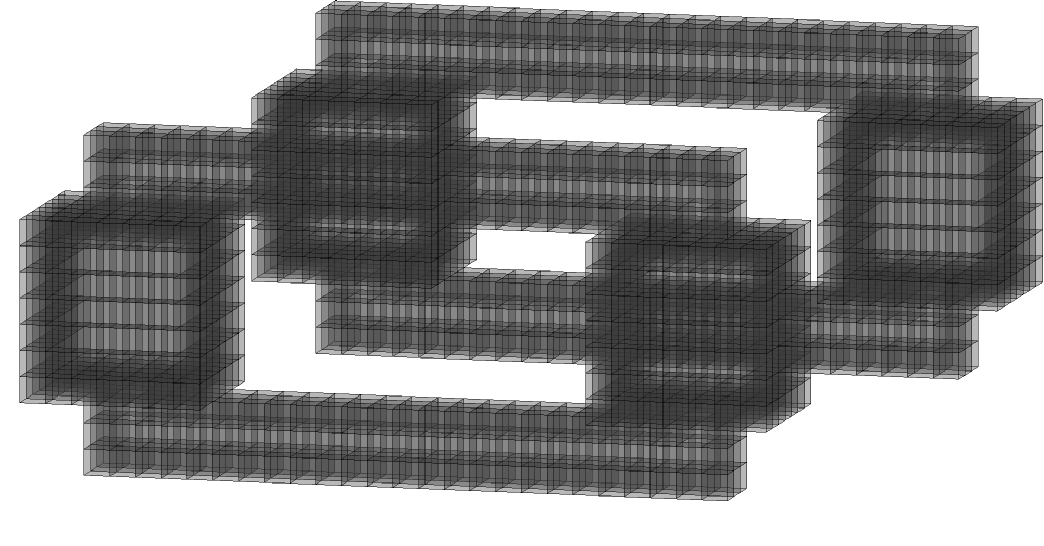}}
\caption{Union of two $3$-dimensional pure cubical complexes, each homotopy equivalent to a torus.  } \label{figsato}
\end{figure}

Figure \ref{figsato} shows a $3$-dimensional pure cubical complex
 $Y$ 
formed from the union of two intersecting pure cubical subcomplexes each of which is 
homotopy equivalent to a torus. The space $Y$ is a union of 1632 cubical $3$-cubes.   
 The four horizontal rectangular tubes of $Y$ have a $3\times 3$ cross section. The four vertical rectangular
tubes of $Y$ have a $7\times 7$ cross section. 
The central axes of the horizontal tubes and vertical tubes lie in a common plane.
 We use $Y$ to construct two $4$-dimensional pure cubical complexes $S$ and $T$ as follows.

\begin{figure}[h]
Profile 1:

\centerline{\includegraphics[height=0.6cm]{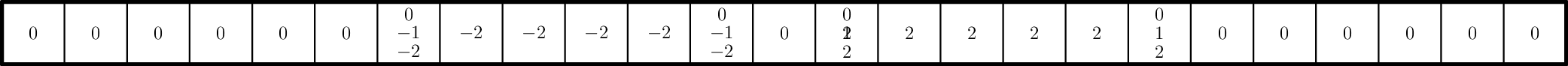}}

\bigskip

Profile 2:

\centerline{\includegraphics[height=0.6cm]{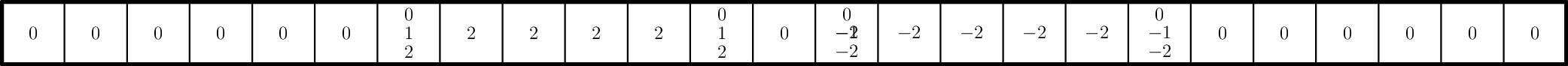}}
\caption{Profiles of cube temperatures for $S$}\label{FIGprofile}
\end{figure}

To construct $S$ we first assign integers $t$, called {\em temperature}, to each $3$-cube of $Y$. To most  $3$-cubes we assign a single temperature $t$, but a
 few  $3$-cubes are assigned three temperatures. All $3$-cubes in the top and bottom horizontal rectangular tubes are assigned one temperature $t=0$. All  $3$-cubes in the four vertical rectangular tubes are again assigned the one temperature $t=0$.
The middle two horizontal rectangular tubes each have length 25. In these two 
tubes the  $x$-coordinate of a cube's  centre determines its temperature(s) 
according to the profiles shown in Figure \ref{FIGprofile}. Profile 1 is the
cross sectional temperature profile of the upper-middle horizontal tube of $Y$;
the lower-middle horizontal tube of $Y$ is assigned Profile 2. The space $S$ is the $4$-dimensional pure cubical complex whose $4$-cubes are centred on the integer vectors $(x,y,z,t)$ where $(x,y,z)$ is the centre of a $3$-cube of $Y$ with temperature $t$.  
 The space $S$ is homotopy equivalent to a disjoint union of two tori.  

\begin{figure}[h]
Profile 1:

\centerline{\includegraphics[height=0.6cm]{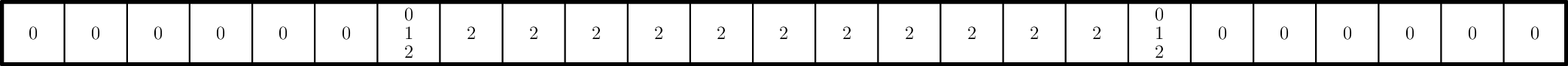}}

\bigskip

Profile 2:

\centerline{\includegraphics[height=0.6cm]{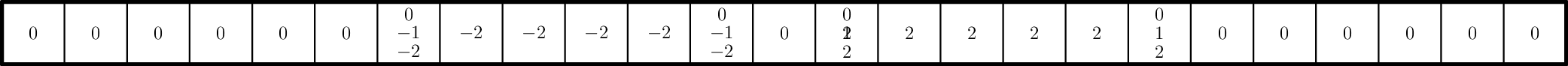}}
\caption{Profiles of cube temperatures for $T$}\label{FIGprofileT}
\end{figure}

The space $T$ is a $4$-dimensional pure cubical complex constructed in the same fashion as $S$ except that the temperature profiles of Figure \ref{FIGprofileT} are used in place of those of Figure \ref{FIGprofile}. The space $T$ is also homotopy equivalent to a disjoint union of two tori.

Let $\mathring{S}$ and $\mathring{T}$ denote the interiors of $S$ and $T$. The complements $\mathbb R^4\setminus \mathring{S}$ and $\mathbb R^4\setminus \mathring{T}$ are subcomplexes of the CW-complex $\mathbb R^4$, both involving infinitely many cells. It is straightforward to construct finite deformation retracts
 $X_S \subset \mathbb R^4\setminus \mathring{S}$ and $X_T
\subset \mathbb R^4\setminus \mathring{T}$ where $X_S$ and $X_T$ are  
$4$-dimensional pure cubical complexes. The computer code of Table \ref{FIGxt}
computes a CW-complex $X_T$ involving 4508573 cells, together with a  sorted list
$Inv(X_T)$  of all possible  abelian invariants of the second
homology groups $H_2(\widetilde X_H,\mathbb Z)$ of 5-fold covering spaces $X_H$ for $X_T$. The list 
$Inv(X_T)$ is an invariant of the homotopy type of $X_T$, and establishes
 that there are $5$-fold covers with 
$H_2(\widetilde X_H,\mathbb Z) = \mathbb Z^{12}$ and $5$-fold
covers with $H_2(\widetilde X_H,\mathbb Z) = \mathbb Z^{16}$.
\begin{table}
\begin{verbatim}
gap> T:=HopfSatohSurface();;
gap> XT:=PureComplexComplement(T);;
gap> XT:=RegularCWComplex(XT);
Regular CW-complex of dimension 4
gap> Size(XT);
4508573

gap> C:=ChainComplexOfUniversalCover(XT);;
gap> L:=Filtered(LowIndexSubgroups(C!.group,5), g->Index(C!.group,g)=5);;
gap> invXT:=List(L,g->Homology(TensorWithIntegersOverSubgroup(C,g),2));;
gap> invXT:=SSortedList(invXT);
[ [ 0, 0, 0, 0, 0, 0, 0, 0, 0, 0, 0, 0 ], 
  [ 0, 0, 0, 0, 0, 0, 0, 0, 0, 0, 0, 0, 0, 0, 0, 0 ] ]
\end{verbatim}
\caption{Computing the CW-complex $X_T$ and an invariant of its homotopy type.} \label{FIGxt}
\end{table}
Similar commands can be applied to $X_S$ to find that for all $5$-fold covers $X_H$ we have $H_2(\widetilde X_H,\mathbb Z) = \mathbb Z^{12}$. Alternatively,
 this information on $X_S$ can be computed more efficiently
using lines 3--7 in the code of Table \ref{TABsmall} which
is based on techniques explained in Section \ref{SECspin}.  Hence $X_T$ is not homotopy equivalent to $X_S$.

Further {\sf GAP} commands establish that $H_0(X_T,\mathbb Z)=\mathbb Z$, 
$H_1(X_T,\mathbb Z)=\mathbb Z^2$,
$H_2(X_T,\mathbb Z)=\mathbb Z^4$,
$H_3(X_T,\mathbb Z)=\mathbb Z^2$,
$H_n(X_T,\mathbb Z)=0$ for $n\ge 4$, and that $\pi_1X_T\cong \pi_1X_S\cong C_\infty\times C_\infty$, $H_n(X_S,\mathbb Z)=H_n(X_T,\mathbb Z)$ for $n\ge 0$.
 These further commands are contained  in the file \texttt{jsc2021-5} detailed
 in  Section \ref{SECrepro}.
So this is an example where the fundamental group and  homology with trivial integral coefficients fail to distinguish between two homotopy inequivalent spaces, but where the second homology with twisted coefficients
does  succeed.

\begin{figure}[h]
\centerline{\includegraphics[height=4cm]{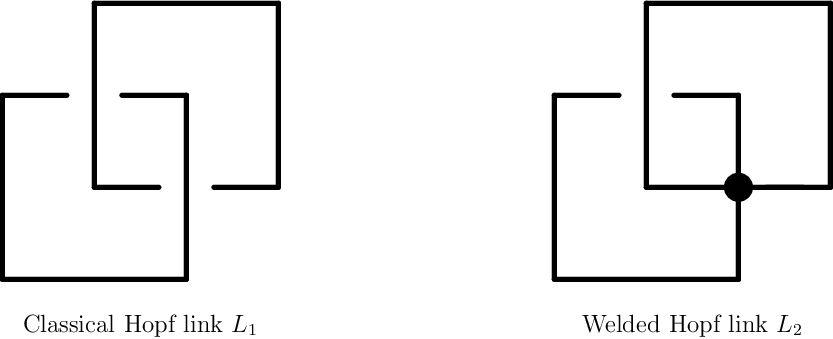}}
\caption{Classical and welded link diagrams} \label{fighopf}
\end{figure}

Figure \ref{fighopf}  shows an example of a classical planar diagram of a link (the Hopf link) and an example of a welded diagram (the welded Hopf link). 
Such a classical or welded link diagram $L$
gives rise to an embedding of a closed surface $tube(L)$ into $\mathbb R^4$ via the {\em tube map}  
of  
\cite{MR1758871}. The surface $tube(L)$ contains one knotted torus $tube(K)$ for each component $K$ of $L$. A good  explanation of the tube map is given on 
Youtube in \cite{davit}. It is shown in \cite{MR1758871} that for a classical 
link diagram $L$ the knotted surface $tube(L)\subset \mathbb R^4$ is the same 
as the surface obtained by spinning ({\em i.e.} rotating)
 the link  about a plane that does not intersect the link. See Section \ref{SECspin} for further details on spun links.  
An algebraic invariant of the homotopy type of the complement $X_L=\mathbb R^4\setminus tube(L)$ is introduced  in \cite{MR2375821}, \cite{MR2506421} and used
in \cite{MR2441256}  to show, for instance, that the classical Hopf link diagram $L_1$ and the welded Hopf link diagram
$L_2$ (see Figure \ref{fighopf}) yield spaces $X_{L_1}$ and
$X_{L_2}$ with different homotopy types.
Consequently the knotted surface $tube(L_1)$ is not ambient isotopic to the knotted surface $tube(L_2)$. Their algebraic invariant is based on J.H.C.\,Whitehead's concept of a crossed module.

The above space $T$ can be viewed as a thickening of two knotted tori in $\mathbb R^4$. The complement $\mathbb R^4\setminus T$ is readily seen to be homeomorphic to the space $X_{L_2}$ associated to the welded Hopf link diagram. 
The complement $\mathbb R^4\setminus S$ is homeomorphic to the space $X_{L_1}$ 
associated to the classical Hopf link diagram. Hence the above computation 
illustrates how the homotopy inequivalence of $X_{L_1}$ and $X_{L_2}$ can be recovered using a {\sf GAP} computation of second homology with local coefficients.
 The computation also yields the following analogue of Theorem 10 of \cite{MR2441256}. For its statement we define the invariant $I_c(\Sigma)$ of a knotted surface $\Sigma\subset \mathbb R^4$ to be 
$$I_c(\Sigma) = \begin{minipage}{11cm}{\rm Set of isomorphism types of abelian groups arising as $H_2(\widetilde X_H,\mathbb Z)$ for some $c$-fold cover
$\widetilde X_H \rightarrow X$ of $X=\mathbb R^4\setminus \Sigma$. }\end{minipage}  $$ 
\begin{theorem}
The invariant $I_c(\Sigma)$ 
of knotted surfaces is powerful enough to distinguish between knotted surfaces
$\Sigma, \Sigma' \subset \mathbb R^4$, with
$\Sigma$
diffeomorphic to 
$\Sigma'$ and 
whose complements have isomorphic fundamental groups and isomorphic integral homology, at least in one particular case.
\end{theorem}

\section{Spinning} \label{SECspin}
As in Subsection \ref{SUBefficient}, let $Y_\ell$ denote a regular CW-complex
 arising from the complement of a link specified by an arc presentation $\ell$. 
The CW-complex $Y_\ell$ is implemented in {\sf HAP}. The first two computer commands in  Table
\ref{TABsmall} construct the complement of the classical Hopf link and 
display  that it involves 303 cells.
\begin{table}
\begin{verbatim}
gap> Y:=KnotComplement([[1,3],[2,4],[1,3],[2,4]]);
Regular CW-complex of dimension 3
gap> Size(Y);
303

gap> SY:=SpunLinkComplement([[1,3],[2,4],[1,3],[2,4]]);;
gap> C:=ChainComplexOfUniversalCover(SY);;
gap> L:=Filtered(LowIndexSubgroupsFpGroup(C!.group,5),g->Index(C!.group,g)=5);;
gap> SSortedList(List(L,g->Homology(TensorWithIntegersOverSubgroup(C,g),2)));
[ [ 0, 0, 0, 0, 0, 0, 0, 0, 0, 0, 0, 0 ] ]

gap> YY:=SimplifiedComplex(Y);
Regular CW-complex of dimension 3
gap> Size(YY);
103
\end{verbatim}
\caption{Computing: (i) a  CW structure on the complement of the Hopf link, (ii) a homotopy invariant of the spun link.} \label{TABsmall}
\end{table}
To explain the remaining commands in Table \ref{TABsmall} we recall details on a spinning construction for links, the origins of which go back to  
 \cite{MR3069446}.
He 
used spinning to construct $4$-dimensional knots from classical knots, but 
we shall give a more general topological description.

Let $X$ be a topological space with subspace $B$. We define the space obtained by {\em spinning $X$ about $B$} to be
$$S_B(X) = X\times \mathbb [0,1] \ / \ \{(x,t)=(x,0) {\rm \ if\ }x\in B {\rm \ ~or\ } t=1\}.$$

An alternative space $\hat S_B(X)$ can be formed from $X$ 
using   
the closed $2$-disk 
 $\mathbb D^2$ and its boundary circle $\mathbb S^1$. We set
$$\hat S_B(X) =  (B\times \mathbb D^2)\  \cup \ (X \times \mathbb S^1) \subset X\times \mathbb D^2.$$ 
 The space $\hat S_B(X)$ is homotopy equivalent to $S_B(X)$. The space
 $\hat S_B(X)$  rather than $ S_B(X)$ is implemented in \textsf{HAP}. More precisely,
suppose that $X$ is a regular CW-complex with subcomplex $B\subset X$. Let $\mathbb D^2$ be given a regular CW structure involving two $0$-cells, two $1$-cells and one $2$-cell. Then the direct product $X\times \mathbb D^2$ is naturally
a regular CW-complex, and $\hat S_B(X)$ is a subcomplex. It is this subcomplex that is implemented.

Artin was interested in the case where $X=\{(x,y,z)\in \mathbb R^3 : z\ge 0\}$ and $B=\{(x,y,z)\in \mathbb R^3 : z= 0\}$. In this case 
$S_B(X)$ is homeomorphic to $\mathbb R^4$, and any knot $\kappa$ embedded in the interior of
$X$ gives rise to a knotted torus $S_\emptyset(\kappa)$ in $\mathbb R^4$. More generally,
 a link $\kappa$ gives rise to an embedded surface $S_\emptyset(\kappa) \subset \mathbb R^4$.
 The complement $\mathbb R^4 \setminus S_\emptyset(\kappa)$ is homeomorphic
to $S_B(X\setminus \kappa)$.

The middle four lines of computer code in Table \ref{TABsmall} show that
for the Hopf link $\kappa$, for $Y=\mathbb R^4\setminus S_\emptyset(\kappa)$, and
for any  $5$-fold covering space $\widetilde Y_H\rightarrow Y$ we have
$H_2(\widetilde Y,\mathbb Z)=\mathbb Z^{12}$. 
 By a result of  \cite{MR1758871} the complement $Y$
of the spun link is homeomorphic to the tube construction for the classical Hopf diagram. 

The final two computer commands in Table \ref{TABsmall} show that Algorithm
\ref{ALGsimplify} can be used to simplify  the CW-structure on the complement  of the classical Hopf link to a CW-structure involving just 103 cells.

\section{The granny and reef knots}\label{SECgranny}
\begin{figure}
\centerline{\includegraphics[height=4cm]{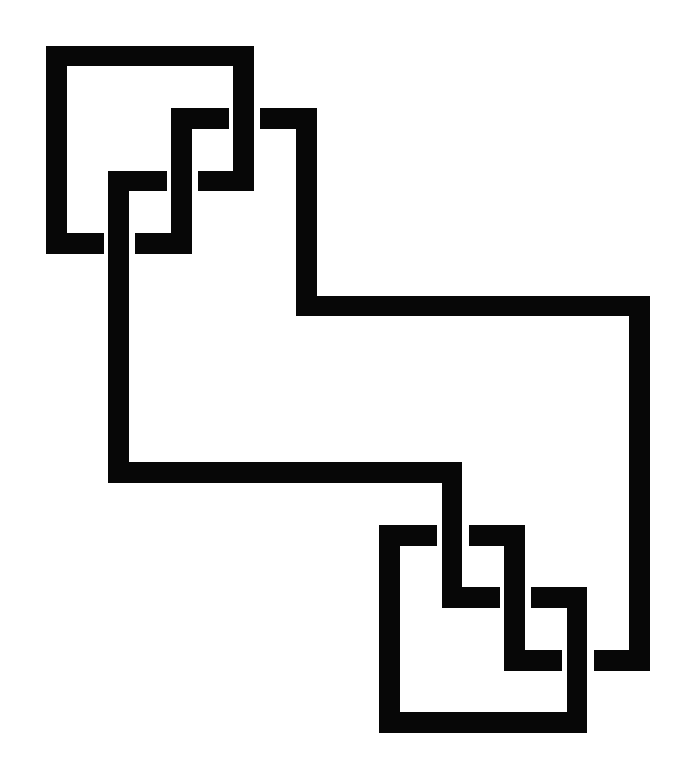}
\includegraphics[height=4cm]{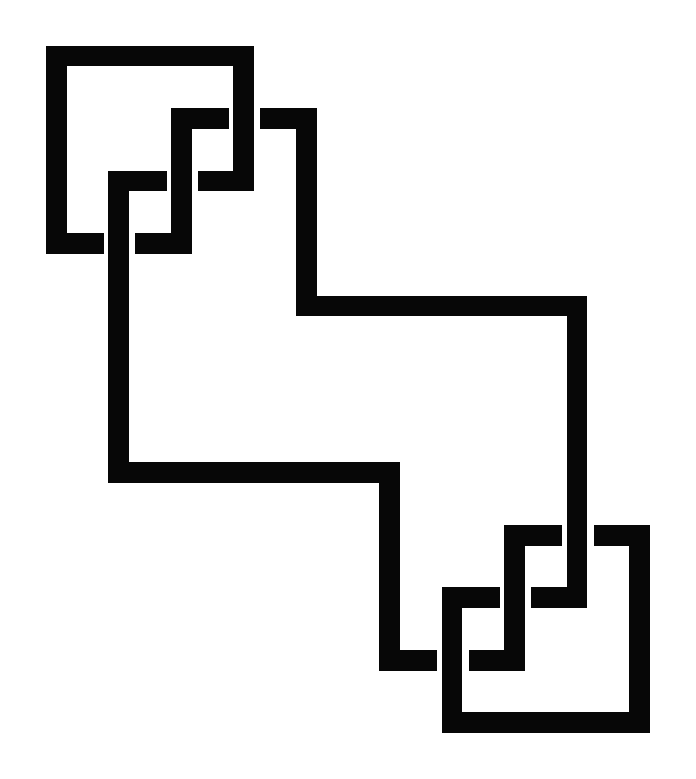}}
\caption{The granny knot (left) and reef knot (right).} \label{FIGgrannyreef}
\end{figure}

Arc diagrams for the granny knot $\kappa_1$ and reef knot $\kappa_2$ 
are shown in Figure \ref{FIGgrannyreef}. Regarding a knot as a $1$-dimensional
subspace of $\mathbb R^3$ we define $X_1=\mathbb R^3\setminus \kappa_1$ and $X_2=\mathbb R^3\setminus \kappa_2$.
Using the method of Subsection \ref{SUBefficient} we can
construct a finite CW-complex $Y_i$ whose interior is homeomorphic to $X_i$, $i=1,2$. The 
\texttt{FundamentalGroup(Y)} function in \textsf{HAP} can be used to compute presentations for $\pi_1Y_i$, $i=1,2$. Tietze moves can be used to establish an isomorphism $\pi_1Y_1 \cong \pi_1Y_2$. Therefore knot invariants, such as the Alexander polynomial, which are based entirely on the knot group are unable to distinguish between the granny and the reef knots. More precisely, they are unable to distinguish between the homeomorphism types of $Y_1$ and $Y_2$.

The integral homology of $Y_1$ and $Y_2$ is readily computed and the spaces turn out to have the same integral homology.

Let us view $X_i$ as the interior of $Y_i$. Let $B'_i$ denote the boundary of 
a 3-dimensional small open tubular neighbourhood of the knot $\kappa_i$. Then $B'_i$ is a 
 2-dimensional CW-subspace of $Y_i$. The subspace $B'_i$ is homeomorphic to a  torus. The \texttt{BoundaryMap(Y)} function in \textsf{HAP} can be used to construct an inclusion $f_i\colon B_i\hookrightarrow Y_i$ of  CW-subspaces
where  $B_i$ is the smallest CW-subspace containing
 those $2$-cells of $Y_i$ that lie in the boundary of exactly one $3$-cell. Thus $B_i$ is a disjoint union of $B'_i$ with a $2$-sphere $\mathbb S^2$. We are interested only in first integral homology $H_1(B_i,\mathbb Z) = H_1(B'_i,\mathbb Z)$ and
so there is no difference if we work  with $B_i$ or with 
$B'_i$,  except that $B_i$ is a little easier to constuct.

For any finite index subgroup $H\le \pi_1Y_i$ the \textsf{HAP} functions
\texttt{U:=UniversalCover(Y)} and
\texttt{p:=EquivariantCWComplexToRegularCWMap(U,H)} can be used
to construct 
the covering map $p_i\colon \widetilde Y_{i,H} \rightarrow Y_i$ which sends
$\pi_1\widetilde Y_{i,H}$ isomorphically onto the subgroup $H$. We can then
use
the function
\texttt{LiftedRegularCWMap(f,p)}  to constuct the subspace
$p_i^{-1}(B_i) \subset 
\widetilde Y_{i,H}$  that maps onto $B_i$, together with the inclusion mapping
$\tilde f_i\colon p_i^{-1}(B_i) \hookrightarrow \widetilde Y_{i,H}$.
 The  induced homology homomorphism
$$\overline f_i\colon H_1(p_i^{-1}(B_i),\mathbb Z) \rightarrow H_1(\widetilde Y_{i,H},\mathbb Z) $$
is readily computed.  This homomorphism $\overline f_i$ is a homeomorphism invariant of $Y_i$. 
In general, the space $p_i^{-1}(B_i)$   is a disjoint union
$p_i^{-1}(B_i) = B_{i,1} \cup B_{i,2} \cup \cdots \cup B_{i,t}$ of path-connected CW-subspaces $B_{i,j}$. The homomorphism $\overline f_i$ is thus a sum of 
homomorphisms 
$$\overline f_{i,j}\colon H_1(B_{i,j},\mathbb Z) \rightarrow H_1(\widetilde Y_{i,H},\mathbb Z) .$$
The collection of all abelian groups ${\rm coker}(\overline f_{i,j})$ arising as the cokernel of $\overline f_{i,j}$ is a homeomorphism invariant of $Y_i$.

For the reef knot $\kappa_2$ and some choice of subgroup $H < \pi_1Y_2$ of index $6$ and some choice of path component $P_{2,j}$, we find that ${\rm coker}(\overline f_{2,j})=\mathbb Z\oplus  \mathbb Z_2 \oplus \mathbb Z_2 \oplus \mathbb Z_8$. However, for the granny knot $\kappa_1$ we find
${\rm coker}(\overline f_{1,j})\not\cong\mathbb Z  \oplus \mathbb Z_2 \oplus \mathbb Z_2 \oplus \mathbb Z_8$ for all subgroups
$H< \pi_1Y_1$ of index $6$ and all path components $P_{1,j}$.  We conclude that $Y_1$ is not homeomorpic to $Y_2$.

These computations motivate the definition of the ambient isotopy invariant
$J_c(\kappa)$ of  a link $\kappa \subset \mathbb R^3$ as 
$$J_c(\kappa) = \begin{minipage}{11cm}{\rm Set of isomorphism types of abelian groups arising as ${\rm coker}(h_j\colon H_1(B_j,\mathbb Z) \rightarrow H_2(\widetilde X_H,\mathbb Z))$ for some $c$-fold cover
$p\colon \widetilde X_H \rightarrow X$ of $X=\mathbb R^3\setminus \kappa$,
and for $B_j$ some path component of $p^{-1}(B)$, where $B$ is the boundary of a small tubular
neigbourhood of $\kappa$. }\end{minipage}  $$
We have proved the following.
\begin{theorem}
The invariant $J_c(\kappa)$
 is powerful enough to distinguish between knots 
$\kappa, \kappa' \subset \mathbb R^3$ 
whose complements have isomorphic fundamental groups and isomorphic integral homology, at least in one particular case.
\end{theorem}

The computations underlying this section can be reproduced by running the example \texttt{jsc2021-6} of Section \ref{SECrepro} in {\sc GAP}.

\section{Tubular neighbourhoods and their complements} \label{SECnbhd}
In this section we consider  a finite regular CW-complex $X$ containing
a  CW-subspace $Y\subset X$, and describe the construction of a finite
CW-complex $W$ that  models
  the complement $X\setminus N_\epsilon(Y)$ of a {\em small open tubular neighbourhood} of $Y$.  We  avoid giving a  precise definition of the open subspace $N_\epsilon(Y) \subset X$ (which is routinely  formulated under the additional assumption that
  $X$ is  piece-wise Euclidean), and focus instead just
on a precise description of $W$. The examples  we have in mind are where $X$ is a contractible compact region of $\mathbb R^n$ for $n=3,4$ and $Y$ is an embedded circle in the case $n=3$ or an embedded surface in the case $n=4$.
 The second author plans to describe applications of a computer implementation of
this construction  in a subsequent paper.

To describe a procedure for constructing the space $W$ we first introduce some terminology and notation for enumerating the cells of $W$ and for describing their homological boundaries.

The CW-complex $W$ will consist of all of the cells in $X\setminus Y$ together with extra cells. We say that a cell of $W$ is {\em internal} if it lies in $X\setminus Y$, and that it is {\em external} otherwise. The complement $X\setminus Y$ is a cell complex -- a union of open cells -- but it is not in general a CW-complex. The external cells ensure that $W$ is a CW-complex.

\begin{figure}
\centerline{\includegraphics[height=3cm]{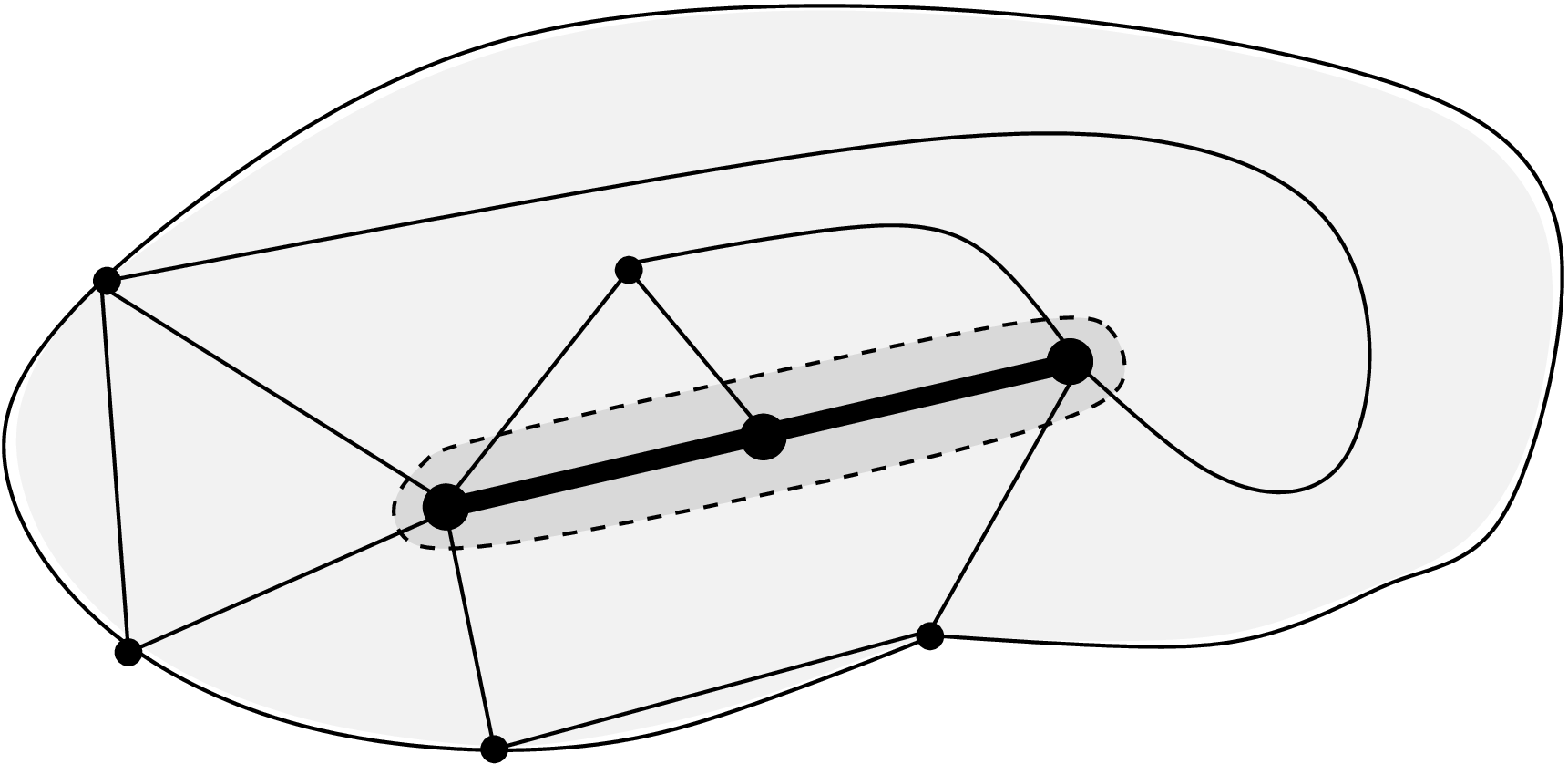}\hspace{1cm} \includegraphics[height=3cm]{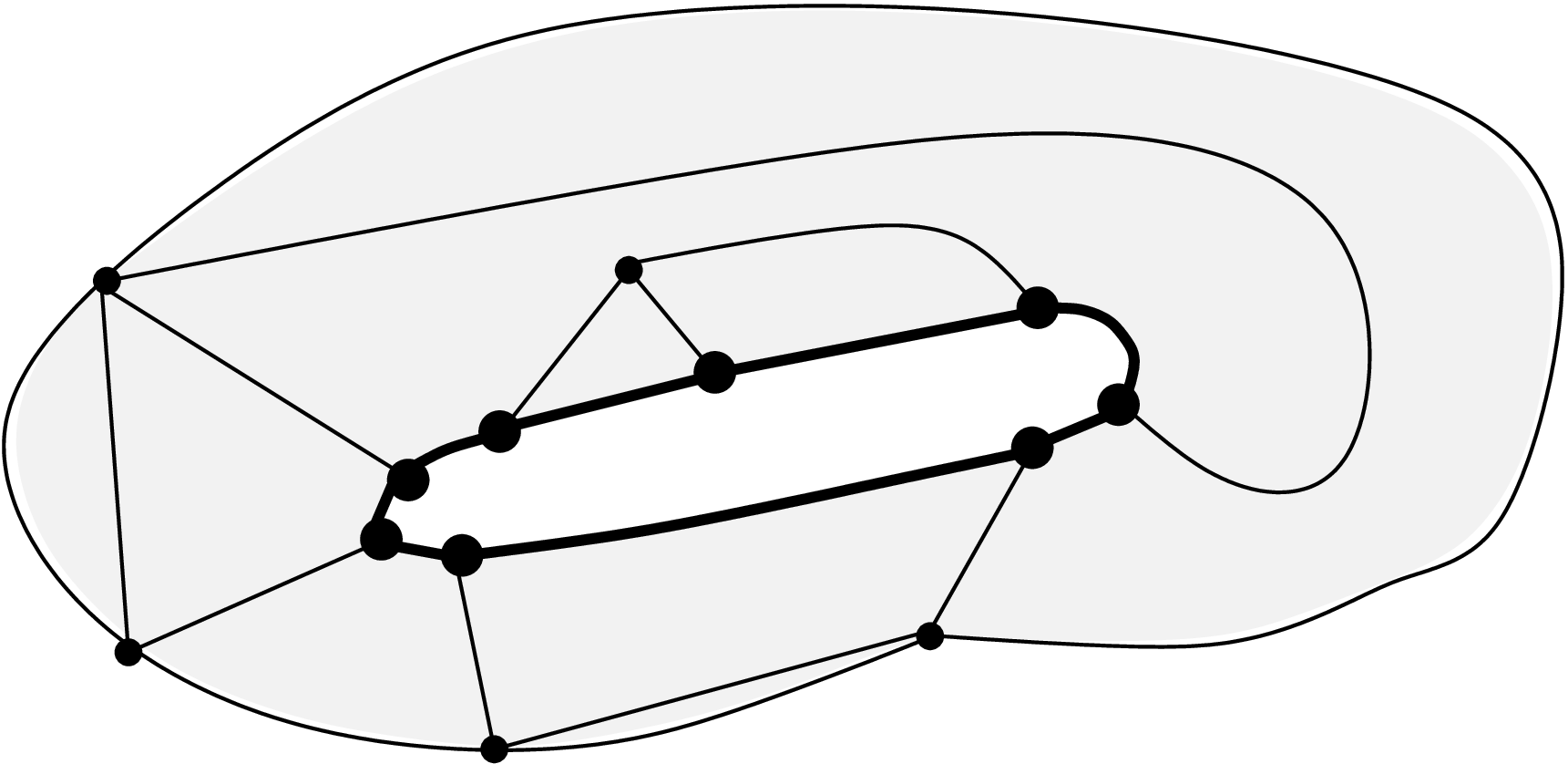}}
\caption{A CW-complex $X$ with open neighbourhood of a subcomplex $Y$ (left). The CW-complex $W$ for this choise of $X$ and $Y$ (right).} \label{FIGnbhd}
\end{figure}

 Figure \ref{FIGnbhd} (left) shows a contractible compact region $X$ in the plane $\mathbb R^2$ endowed with a CW-structure involving $8$ $2$-cells and various cells of dimension $1$ and $0$. A $1$-dimensional  subcomplex $Y \subset X$, involving two $1$-cells and three $0$-cells,
is shown in bold. An open tubular neighbourhood  $N_\epsilon(Y) \subset X$ is shown in dark grey.
Figure \ref{FIGnbhd} (right) shows the CW-complex $W$ corresponding to
 this particular choice of $X$ and $Y$.

The  closure $\overline{e^n}$ of any $n$-cell $e^n \subset X$ is
a  finite CW-subcomplex of $X$.
For an internal $n$-cell  $e^n\subset X\setminus Y$ the intersection
$\overline e^n \cap Y$ is also a CW-subcomplex of $X$, which we can express as a union
$$\overline e^n \cap Y = A_1^{e^n} \cup A_2^{e^n} \cup \cdots A_k^{e^n}$$
of its path-components $A_i^{e^n}$. Each $A_i^{e^n}$ is a CW-subcomplex of $Y$.

Let us make the following simplifying assumption. When a space fails to meet the assumption it is sometimes possible to apply  barycentric subdivision to offending top-dimensional cells (or possibly some less costly subdivision) in order to modify the  CW structure so that it does satisfy the assumption.

{\bf Simplifying assumption:} Let us suppose that each $A_i^{e^n}$ is contractible.

Under this assumption the CW-complex $W$ has precisely one internal $n$-cell $e^n$ for each cell
$e^n\subset X\setminus Y$, and  precisely one external $(n-1)$-cell
$$f^{e^n}_{A_i}$$
for each path-component $A_i$ of $\overline e^n \cap Y$.

The cellular chain complex $C_\ast W$ is generated in degree $n$ by the $n$-cells $e^n$ of $W$. By the {\em homological boundary} of $e^n$ we mean
 its image $de^n$ under the boundary homomorphism $d\colon C_nW \rightarrow C_{n-1}W$. The homological boundary $de^n$  is a sum of finitely many distinct
 $(n-1)$-cells each occuring with coefficient $\pm 1$. We'll abuse terminology and say that $de^n$ is a sum of these distinct  $(n-1)$ cells (without mentioning that  coefficients are $\pm 1$ in the sum).

The homological boundary of an internal cell $e^n$ in $W$ is  a sum of all the internal $(n-1)$-cells lying in $\overline e^n$ together with all external $(n-1)$-cells $f^{e^n}_{A_i}$ with $A_i$ a path-components of $\overline e^n \cap Y$.
The homological boundary of an external $(n-1)$-cell $f^{e^n}_{A_i}$ is a sum of all those external $(n-1)$-cells $f^{e^{n-1}}_B$ with $e^{n-1}$ an $(n-1)$-cell of $\overline e^n$ and path-component $B$ of $\overline e^n\cap Y$ with $B\subseteq A_i$.

\begin{algorithm}
\caption{Excise an open tubular neighbourhood from a regular CW-complex.}
\label{ALGnbhd}
\begin{algorithmic}[1]
\Require A finite regular  CW-complex $X$ and subcomplex $Y\subset X$ satisfying the above simplifying assumption.
\Ensure A  finite regular CW-complex $W$.
\Procedure{}{}

\State For each internal cell $e^n\subset X\setminus Y$ compute the CW-complex
$\overline e^n \cap Y$ and express it as a union $e^n \cap Y = A_1^{e^n} \cup A_2^{e^n} \cup \cdots A_k^{e^n}$
of its path-components $A_i^{e^n}$. This determines the number of cells of each dimension in $W$.

\State Create a list of empty lists $B=[[ ], [ ], \cdots, [ ]]$ of length
$1+\dim X$.

\State  For $0\le n\le \dim X$ set $B[n+1]=[b_1, b_2, \cdots, b_{\alpha_n}]$ where $\alpha_n$ denotes the number of $n$-cells in $W$ and $b_i$ is a list of integers determining the $(n-1)$-cells with non-zero coefficient in the homological
boundary of the $i$-th $n$-cell of $W$. In the ordering of $n$-cells it is convenient to order all internal cells before any external cell.

\State Represent the data in $B$ as a regular CW-complex $W$ and return $W$.
\EndProcedure
\end{algorithmic}
\end{algorithm}

\section{A homotopy classification of maps}\label{SECclass}
Our original motivation for attempting computations of local cohomology is 
a classical result of Whitehead
concerning the set $[W,X]_\phi$ of based homotopy classes of maps
$f\colon W\rightarrow X$ between connected
CW-complexes $W$ and $X$  that induce a
given homomorphism $\phi\colon \pi_1W\rightarrow \pi_1X$ of fundamental groups.
As mentioned in the Introduction, we assume that the spaces $W$, $X$  have preferred base points, that the
maps $f$ preserve base points, and that base points are preserved at each stage of  a homotopy between maps.
 The homotopy groups of a space $X$ are denoted by $\pi_nX$, $n\ge 0$.
A chain map $\chi_\ast\colon C_\ast\widetilde W \rightarrow C_\ast\widetilde X$
between cellular chain complexes of universal covers 
is said to be {\em $\phi$-equivariant} if $\chi_n(gw) = \phi(g)\chi_n(w)$
for all $g\in \pi_1 W$, $w\in C_n\widetilde W$, $n\ge 0$. 
Two such $\phi$-equivariant chain maps $\chi_\ast, \psi_\ast
\colon C_\ast\widetilde W \rightarrow C_\ast\widetilde X$ are said to be 
{\em based chain homotopic}
if
there exist  $\phi$-equivariant homomorphisms $h_i\colon C_i \widetilde W \rightarrow C_{i+1}\widetilde X$ for $i\ge 0$,
 such that
\begin{equation}
 \chi_i(w) - \psi_i(w) = d_{i+1} h_i(w) + h_{i-1}d_i(w)
\end{equation}
for $w\in C_i\widetilde W$ and $h_{-1}=0$. 
Based chain homotopy is an equivalence relation and we denote by
$[C_\ast\widetilde W,C_\ast\widetilde X]_\phi$ the set of based chain homotopy classes of $\phi$-equivariant chain maps 
$C_\ast\widetilde W \rightarrow C_\ast\widetilde X$.
The following theorem is proved in \cite{MR0030760}.

\begin{theorem}
[Whitehead] \label{thm4}
Let $W,X$ be connected CW-complexes, and let $\phi\colon \pi_1W \rightarrow \pi_1X$ be a group homomorphism. Suppose that $W$ is of dimension $n$ and that $\pi_iX=0$ for $2\le i\le n-1$. Then there is a bijection
$$[W,X]_\phi \cong [C_\ast\widetilde W,C_\ast\widetilde X]_\phi .$$
\end{theorem}

To derive Theorem \ref{thm} from Theorem \ref{thm4} we need to establish a bijection
$[C_\ast\widetilde W,C_\ast\widetilde X]_\phi \cong H^n(W, H_n(\widetilde X))$.
Let 
  $W,X$ and $\phi$ satisfy the hypothesis of Theorem \ref{thm}.
Without loss of generality, we impose the additional hypothesis that both $W$ and $X$ have precisely one   $0$-cell, say $e^0\in W$, $e'^0\in X$, and that these $0$-cells are taken to be the base points. 
 Furthermore, we identify cells of $W, X$ with free $\mathbb Z\pi_1$-generators of $C_\ast\widetilde X, C_\ast\widetilde W$. 
 Let $f\colon W \rightarrow X$ be a cellular map inducing $\phi$. Such a map
$f$ exists because, by the theory of group presentations, the  $2$-skeleta $W^2$, $X^2$ correspond to group presentations of $\pi_1W$, $\pi_1X$;  the homomorphism $\phi$ is induced by at least one cellular map $f\colon W^2 \rightarrow X^2$. This map of 2-skeleta extends to a cellular map $f\colon W \rightarrow X$
since $\pi_iX=0$ for $2\le i\le n-1$. 
 Let $\chi_\ast\colon C_\ast\widetilde W \rightarrow C_\ast\widetilde X$ be the chain map induced by $f$. 
  Note that $\chi_\ast$ is pointed in the sense that $\chi_0(e^0) =  e'^0$. 
Let $[\chi_\ast]$ denote the based homotopy class of
$\chi_\ast$. Then $[\chi_\ast] \in 
[C_\ast\widetilde W,C_\ast\widetilde X]_\phi$. We choose 
(one such) $\chi_\ast$ as the preferred $\phi$-equivariant chain map.

Suppose that $\psi_\ast\colon C_\ast\widetilde W \rightarrow C_\ast\widetilde X$ is some other $\phi$-equivariant pointed chain map.
Then $\chi_0=\psi_0$. Suppose that for some $0\le k\le n-2$ we have 
$\phi_i=\psi_i$ for $0\le i\le k$. Then for each free $\mathbb Z\pi_1W$-generator 
$e^k_j \in  C_k\widetilde W$ we have $\chi_k(e^k_j) - \psi_k(e^k_j) \in
\ker(d_k\colon C_k\widetilde X \rightarrow C_{k-1}\widetilde X$).
 Since $H_k(C_\ast\widetilde X)=\pi_k(X)=0$ we can choose some element
$h(e'^{k}_j) \in C_{k+1}\widetilde X$ satisfying $\chi_k(e^k_j) - \psi_k(e^k_j) =d_{k+1}h(e'^{k}_j)$. Here $h(e'^{k}_j)$ is just our label for the chosen element.
 This defines a $\phi$-equivariant homomorphism
$h\colon C_k\widetilde W \rightarrow C_{k+1}\widetilde X$.

We define a $\phi$-equivariant chain map $\psi'_\ast\colon C_\ast\widetilde W \rightarrow C_\ast \widetilde X$ by setting 
$\psi'_i=\psi_i$ $ (i\ne k,k+1)$,
$\psi'_k=\chi_k$, $\psi'_{k+1}= \psi_{k+1} + hd_{k+1}$. 
Then $\psi'_\ast$ is 
based chain homotopic to $\psi_\ast$ and $\psi'_i=\chi_i$ for $0\le i\le k+1$.
By induction the based homotopy class $[\psi_\ast]$ is represented by a chain map $\psi_\ast$ satisfying $\psi_i=\chi_i$ for $0\le i\le n-1$.
 Using this representative we define the $\phi$-equivariant
cocycle
$$C_{[\chi_\ast]-[\psi_\ast]} = \chi_n - \psi_n \colon C_n\widetilde W \rightarrow H_n(C_\ast \widetilde X) .$$
For fixed $\chi_\ast$, the assignment $\psi_\ast \mapsto C_{[\chi_\ast]-[\psi_\ast]} $ induces the  bijection $[C_\ast\widetilde W,C_\ast\widetilde X]_\phi \cong H^n(W, H_n(\widetilde X))$ which proves Theorem
\ref{thm}.

\section{Reproducibility of computations}\label{SECrepro}
All computations in the paper were obtained using {\sc GAP} v4.11.0 and {\sc HAP} v1.29  and  Debian GNU/Linux 10 (buster) on a
laptop with an  i7--6500U CPU and 4GB of RAM.
The code in Tables 1--4 has been stored as examples in {\sc HAP} which
 can be run using the  command

\begin{verbatim}
gap> HapExample("jsc2021-n"); 
\end{verbatim}

\noindent for $n=1,2,3,4$. Two further examples relating to computations in the paper have also been stored in \texttt{HAP}. Timings for the  examples, in
\texttt{hour:min:sec} format, are as follows.

\bigskip
\begin{tabular}{l|p{9cm}|l}
 Example & Details & Time  \\
\hline \texttt{jsc2021-1} & Code of Table 1.& 0:06:09\\
\texttt{jsc2021-2} & Code of Table 2. Requires \texttt{jsc2021-1}. &0:43:13\\
\texttt{jsc2021-3} & Code of Table 3.& 0:18:36\\
\texttt{jsc2021-4} & Code of Table 4.& 0:00:02\\
\texttt{jsc2021-5} & Establish that $X_S$, $X_T$ have identical integral homology and fundamental group. Requires \texttt{jsc2021-3} \& \texttt{jsc2021-4}.&0:00:06\\
\texttt{jsc2021-6} & Computations underlying Section 9. &0:04:08
\end{tabular}

\section{Acknowledgments} The second author is grateful to the Irish Research Council for a postgraduate scholarship under project GOIPG/2018/2152.
 The authors thank the referee for helpful comments and for suggesting the inclusion of Section 12.

\bibliographystyle{elsarticle-harv}

\bibliography{mybib}{}

\def\cprime{$'$}
\begin{thebibliography}{30}
\expandafter\ifx\csname natexlab\endcsname\relax\def\natexlab#1{#1}\fi
\expandafter\ifx\csname url\endcsname\relax
  \def\url#1{\texttt{#1}}\fi
\expandafter\ifx\csname urlprefix\endcsname\relax\def\urlprefix{URL }\fi

\bibitem[{Artin(1925)}]{MR3069446}
Artin, E., 1925. Zur {I}sotopie zweidimensionaler {F}l\"{a}chen im {$R_4$}.
  Abh. Math. Sem. Univ. Hamburg 4~(1), 174--177.
\newline\urlprefix\url{https://doi-org.libgate.library.nuigalway.ie/10.1007/BF02950724}

\bibitem[{Brendel et~al.(2015)Brendel, D{\l}otko, Ellis, Juda, and
  Mrozek}]{EllisMrozek}
Brendel, P., D{\l}otko, P., Ellis, G., Juda, M., Mrozek, M., 2015. Computing
  fundamental groups from point clouds. Appl. Algebra Engrg. Comm. Comput.
  26~(1-2), 27--48.
\newline\urlprefix\url{http://dx.doi.org/10.1007/s00200-014-0244-1}

\bibitem[{Davit(2017)}]{davit}
Davit, E., 2017. Knots in 4d -- part 3. Youtube.
\newline\urlprefix\url{https://www.youtube.com/watch?v=p5tPb5yN-CE}

\bibitem[{Ellis(2019)}]{ellisbook}
Ellis, G., 2019. An invitation to computational homotopy. Oxford University
  Press, UK.
\newline\urlprefix\url{https://global.oup.com/academic/product/an-invitation-to-computational-homotopy-9780198832980}

\bibitem[{Ellis(November 2019)}]{hap}
Ellis, G., November 2019. HAP -- Homological Algebra Programming, Version 1.21.
  (\texttt{http://www.gap-system.org/Packages/hap.html}).

\bibitem[{Ellis and Fragnaud(2018)}]{MR3896314}
Ellis, G., Fragnaud, C., 2018. Computing with knot quandles. J. Knot Theory
  Ramifications 27~(14), 1850074, 18.
\newline\urlprefix\url{https://doi-org.libgate.library.nuigalway.ie/10.1142/S0218216518500748}

\bibitem[{Ellis and Hegarty(2014)}]{ELLISHEGARTY}
Ellis, G., Hegarty, F., 2014. Computational homotopy of finite regular
  {CW}-spaces. J. Homotopy Relat. Struct. 9~(1), 25--54.
\newline\urlprefix\url{http://dx.doi.org/10.1007/s40062-013-0029-4}

\bibitem[{Forman(1998)}]{MR1612391}
Forman, R., 1998. Morse theory for cell complexes. Adv. Math. 134~(1), 90--145.
\newline\urlprefix\url{http://dx.doi.org/10.1006/aima.1997.1650}

\bibitem[{Forman(2002)}]{MR1939695}
Forman, R., 2002. A user's guide to discrete {M}orse theory. S\'em. Lothar.
  Combin. 48, Art.\ B48c, 35.

\bibitem[{Fox(1954)}]{MR62125}
Fox, R.~H., 1954. Free differential calculus. {II}. {T}he isomorphism problem
  of groups. Ann. of Math. (2) 59, 196--210.
\newline\urlprefix\url{https://doi-org.libgate.library.nuigalway.ie/10.2307/1969686}

\bibitem[{GAP(2013)}]{GAP4}
GAP, 2013. {GAP -- Groups, Algorithms, and Programming, Version 4.5.6}. The
  GAP~Group, (\texttt{http://www.gap-system.org}).

\bibitem[{Hatcher(2002)}]{opacb1122188}
Hatcher, A., 2002. Algebraic topology. Cambridge University Press, Cambridge,
  New York, autre(s) tirage(s) : 2003,2004,2005,2006.
\newline\urlprefix\url{http://opac.inria.fr/record=b1122188}

\bibitem[{Hempel(1984)}]{MR739142}
Hempel, J., 1984. Homology of coverings. Pacific J. Math. 112~(1), 83--113.
\newline\urlprefix\url{http://projecteuclid.org.libgate.library.nuigalway.ie/euclid.pjm/1102710101}

\bibitem[{Jones(1988)}]{MR968920}
Jones, D., 1988. A general theory of polyhedral sets and the corresponding
  {$T$}-complexes. Dissertationes Math. (Rozprawy Mat.) 266, 110.

\bibitem[{Kauffman and Faria~Martins(2008)}]{MR2441256}
Kauffman, L.~H., Faria~Martins, J.~a., 2008. Invariants of welded virtual knots
  via crossed module invariants of knotted surfaces. Compos. Math. 144~(4),
  1046--1080.
\newline\urlprefix\url{https://doi.org/10.1112/S0010437X07003429}

\bibitem[{Kr\u{a}\'{z}\'{a}l et~al.(2013)Kr\u{a}\'{z}\'{a}l, Matou\v{s}ek, and
  Sergeraert}]{Krazal:2013}
Kr\u{a}\'{z}\'{a}l, M., Matou\v{s}ek, J., Sergeraert, F., Dec. 2013.
  Polynomial-time homology for simplicial {E}ilenberg--{M}aclane spaces. Found.
  Comput. Math. 13~(6), 935--963.
\newline\urlprefix\url{https://doi.org/10.1007/s10208-013-9159-7}

\bibitem[{Martins(2007)}]{MR2375821}
Martins, J.~a.~F., 2007. Categorical groups, knots and knotted surfaces. J.
  Knot Theory Ramifications 16~(9), 1181--1217.
\newline\urlprefix\url{https://doi-org.libgate.library.nuigalway.ie/10.1142/S0218216507005713}

\bibitem[{Martins(2009)}]{MR2506421}
Martins, J.~F., 2009. The fundamental crossed module of the complement of a
  knotted surface. Trans. Amer. Math. Soc. 361~(9), 4593--4630.
\newline\urlprefix\url{http://dx.doi.org/10.1090/S0002-9947-09-04576-0}

\bibitem[{Mazur(1959)}]{MR117693}
Mazur, B., 1959. On embeddings of spheres. Bull. Amer. Math. Soc. 65, 59--65.
\newline\urlprefix\url{https://doi.org/10.1090/S0002-9904-1959-10274-3}

\bibitem[{Ocken(1990)}]{MR984809}
Ocken, S., 1990. Homology of branched cyclic covers of knots. Proc. Amer. Math.
  Soc. 110~(4), 1063--1067.
\newline\urlprefix\url{https://doi-org.libgate.library.nuigalway.ie/10.2307/2047757}

\bibitem[{Rees and Soicher(2000)}]{MR1743389}
Rees, S., Soicher, L., 2000. An algorithmic approach to fundamental groups and
  covers of combinatorial cell complexes. J. Symbolic Comput. 29~(1), 59--77.
\newline\urlprefix\url{http://dx.doi.org/10.1006/jsco.1999.0292}

\bibitem[{Satoh(2000)}]{MR1758871}
Satoh, S., 2000. Virtual knot presentation of ribbon torus-knots. J. Knot
  Theory Ramifications 9~(4), 531--542.
\newline\urlprefix\url{https://doi.org/10.1142/S0218216500000293}

\bibitem[{Spanier(1981)}]{spanier}
Spanier, E., 1981. Algebraic Topology, corrected reprint of the 1966 original.
  Springer, New York.

\bibitem[{{The Sage Developers}(2019)}]{palmieri}
{The Sage Developers}, 2019. {S}ageMath, the {S}age {M}athematics {S}oftware
  {S}ystem ({V}ersion 8.7.0). {\tt
  http://doc.sagemath.org/html/en/reference/homology/sage/homology/simplicial\_set.html}.

\bibitem[{Trotter(1962)}]{MR143201}
Trotter, H.~F., 1962. Homology of group systems with applications to knot
  theory. Ann. of Math. (2) 76, 464--498.
\newline\urlprefix\url{https://doi.org/10.2307/1970369}

\bibitem[{\v{C}adek et~al.(2014{\natexlab{a}})\v{C}adek, Kr\v{c}\'{a}l,
  Matou\v{s}ek, Sergeraert, Vok\v{r}\'{\i}nek, and Wagner}]{CADEK:2014}
\v{C}adek, M., Kr\v{c}\'{a}l, M., Matou\v{s}ek, J., Sergeraert, F.,
  Vok\v{r}\'{\i}nek, L., Wagner, U., Jun. 2014{\natexlab{a}}. Computing all
  maps into a sphere. J. ACM 61~(3), 17:1--17:44.
\newline\urlprefix\url{http://doi.acm.org/10.1145/2597629}

\bibitem[{\v{C}adek et~al.(2014{\natexlab{b}})\v{C}adek, Kr\v{c}\'{a}l,
  Matou\v{s}ek, Vok\v{r}\'{\i}nek, and Wagner}]{kcral3}
\v{C}adek, M., Kr\v{c}\'{a}l, M., Matou\v{s}ek, J., Vok\v{r}\'{\i}nek, L.,
  Wagner, U., 2014{\natexlab{b}}. Polynomial-time computation of homotopy
  groups and {P}ostnikov systems in fixed dimension. SIAM Journal on Computing
  43~(5), 1728--1780.

\bibitem[{Whitehead(1949)}]{MR0030760}
Whitehead, J. H.~C., 1949. Combinatorial homotopy. {II}. Bull. Amer. Math. Soc.
  55, 453--496.

\bibitem[{Whitehead(1950)}]{MR0035437}
Whitehead, J. H.~C., 1950. Simple homotopy types. Amer. J. Math. 72, 1--57.

\bibitem[{Zeeman(1961)}]{MR123335}
Zeeman, E.~C., 1961. Knotting manifolds. Bull. Amer. Math. Soc. 67, 117--119.
\newline\urlprefix\url{https://doi.org/10.1090/S0002-9904-1961-10529-6}

\end{thebibliography}

\end{document}